\documentclass[review]{elsarticle}

\usepackage{lineno,hyperref}
\modulolinenumbers[5]
\usepackage{amsmath}
\graphicspath{{Images1/}}
\usepackage{caption}
\usepackage{subcaption}
\usepackage{epstopdf}
\usepackage[margin=2.5cm]{geometry}
\usepackage{lipsum}
\journal{Mathematical Biosciences}
\usepackage{comment}
\usepackage{floatrow}
\usepackage[labelformat=empty, position=top]{subcaption}
\usepackage[export]{adjustbox}

\begin{document}

\begin{frontmatter}

\title{The role of viral infectivity in oncolytic virotherapy outcomes: a mathematical study}

\author[usydmaths]{Pantea Pooladvand}
\author[Hanyang,Nano]{Chae-Ok Yun}
\author[Hanyang,Nano]{A-Rum Yoon}
\author[usydmaths]{Peter S. Kim}
\author[swinburne]{Federico Frascoli}
\cortext[cor]{Pantea Pooladvand}
\address[usydmaths]{School of Mathematics and Statistics, The University of Sydney, Sydney, NSW 2006, Australia}
\address[Hanyang]{Department of Bioengineering, Collage of Engineering, Hanyang University, Seoul, South Korea}
\address[swinburne]{Department of Mathematics, Faculty of Science, Engineering and Technology, Swinburne University of Technology, Melbourne, VIC 3122, Australia}
\address[Nano]{Institute of Nano Science and Technology (INST), Hanyang University, Seoul, South Korea}


\begin{abstract}
A model capturing the dynamics between virus and tumour cells in the context of oncolytic virotherapy is presented and analysed. The ability of the virus to be internalised by uninfected cells is described by an infectivity parameter, which is inferred from available experimental data. The parameter is also able to describe the effects of changes in the tumour environment that affect viral uptake from tumour cells.

Results show that when a virus is inoculated inside a growing tumour, strategies for enhancing infectivity do not lead to a complete eradication of the tumour. Within typical times of experiments and treatments, we observe the onset of oscillations, which always prevent a full destruction of the tumour mass. These findings are in good agreement with available laboratory results.

Further analysis shows why a fully successful therapy cannot exist for the proposed model and that care must be taken when designing and engineering viral vectors with enhanced features. In particular, bifurcation analysis reveals that creating longer lasting virus particles or using strategies for reducing infected cell lifespan can cause unexpected and unwanted surges in the overall tumour load over time. Our findings suggest that virotherapy alone seems unlikely to be effective in clinical settings unless adjuvant strategies are included.
\end{abstract}


\end{frontmatter}


\section{Introduction}

An oncolytic, or anti-tumour virus is a type of virus that preferentially or exclusively targets, infects and kills tumour cells. Once infected cells burst (a process known as lysis), a new population of viruses, resulting from virus replication inside the infected cells, is released into the tumour environment. Since the first experiment approximately seventy years ago \cite{kelly2007history}, a number of fundamental and technological advances have led to a host of oncolytic viruses currently used in human trials and animal models \cite{patel2013oncolytic} with the goal of creating viral vectors that are able to eradicate growing tumours. Examples include reovirus, measles, herpes simplex virus (HSV) and adenovirus \cite{eissa2018current}.

Despite some partial success, implementation of virotherapy as a routine method for treating solid tumours is not in sight. Limitations in treatment include inconsistent expression levels of receptors on tumour cells to allow virus binding and infection, such as the coxsackie Ad receptor (CAR) essential for the internalisation of adenovirus \cite{kim2012enhancing}; viral clearance by immune cells \cite{alvarez2012nk}; physical barriers such as interstitial fluid pressure \cite{stohrer2000oncotic}; and the extracellular matrix (EMC), a key element in the tumour stroma. \cite{jain2010delivering}. These limitations result in reduced infectivity and lack of efficacy of virotherapy in many types of solid tumours.

Experimental studies have utilised an array of methods to increase viral infectivity in solid tumours. McKee {\em et al.} show through multiphoton imaging that viral distribution is impeded by collagen and that co-injection of collagenase results in approximately a three-fold increase in the area of viral distribution than with viral particles alone \cite{mckee2006degradation}. Kim {\em et al.} also compare the efficacy of treatment between standard oncolytic adenoviruses and adenoviruses that express relaxin, which is a protein that degrades collagen \cite{kim2006relaxin}. They found that relaxin-expressing viruses consistently outperform the ordinary ones by infecting considerably larger areas of collagen-dense tumours. Similar results relating the degradation of collagen to increased efficacy of viral treatment and delay of tumour growth have been produced using co-delivery of other forms of enzymes   \cite{ganesh2008intratumoral,guedan2010hyaluronidase,diop2011losartan}. More recently, strategies targeting cancer-associated fibroblasts, responsible for synthesis of collagen in the tumour environment demonstrated greater infiltration of viruses and tumour suppression \cite{yu2017t, de2019targeting}. To improve infectivity and viral uptake, adenoviruses are being genetically modified to overcome CAR-dependent internalisation~\cite{yoon2015vesicular}. Another promising approach is the combination of oncolytic virus with chemotherapy drugs known as chemovirotherapy. Gomez-Gutierrez {\em et al.} \cite{gomez2016combined} showed that the combined effect of chemotherapeutic drug temozolomide (TMZ) and oncolytic adenovirus (adeAdhz60) has a synergistic killing effect on three lung cancer cell lines. Combination therapies also include  immunosuppressive drugs, administered in conjunction with virotherapy to reduce viral clearance by immune cells \cite{meisen2015impact}. On the other hand, the anti-viral response from the immune system during oncolytic virotherapy has also been utilised to activate the immune system against tumour cells~\cite{russell2019oncolytic, gujar2018antitumor}. Reviews of oncolytic virus limitations and advancements can be found here \cite{zheng2019oncolytic, de2020oncolytic, hong2019overcoming}.

Oncolytic virotherapy has been recently investigated in the biomathematical literature with the aid of different approaches. Models exploring viral infectivity and spread in tumours have described lack of diffusion, oscillatory behaviour and conditions for tumour eradication using ordinary differential equations (ODEs)~\cite{phan2017role, wodarz2004computational, wodarz2009towards, tian2011replicability}, partial differential equations (PDEs) \cite{wu2001modeling, wein2003validation, Fried2006}, agent-based~\cite{rodriguez2017complex, jenner2020enhancing} and hybrid models. Mok {\em et al.} \cite{mok2009mathematical}, uses a system of PDEs to track free, bound and internalised viruses to describe the spread of herpes simplex virus when administered by injection to a solid tumour. The model explores the role of virus binding, internalisation, diffusion and degradation. The authors conclude that rapid binding, internalisation and lack of diffusion due to tumour density impede virus propagation. The model highlights that modifications such as alterations in the viral envelope to decrease binding affinity and degradation of the extracellular matrix to increase diffusivity could improve treatment outcome. Wodarz {\em et al.} \cite{wodarz2012complex} use an agent-based model to recreate spatial patterns produced in vitro by viral infection, when human embryonic kidney cells are infected with oncolytic adenovirus. The model is then used to predict the long-term treatment outcome based on the spatial patterns. The authors conclude that one of the spatial patterns leads to tumour extinction while two others fail. Malinzi {\em et al.} \cite{malinzi2017modelling} use a system of PDEs to study the spatiotemporal distribution of viruses in  chemovirotherapy. The model includes uninfected and infected tumour cells, free viruses and chemotherapy drug. Parameter sensitivity analysis is performed to determine the key drivers of cancer remission during treatment. The authors also study the temporal model which suggests that virus burst size and infection rate determine virotherapy outcome. Friedman and Lai \cite{friedman2018combination} develop a PDE model to study the combined therapy of oncolytic viruses and immune checkpoint inhibitors. Immune checkpoints stop T cells from attacking cancer cells and inhibiting these checkpoints has clinically shown greater efficacy in treatment~\cite{rajani2015harnessing}. In this model, Friedman and Lai consider populations of uninfected cancer cells, infected cancer cells, extracellular and intracellular viruses, macrophages, dendritic and CD4+ T cells. The authors explore the efficacy of the combined treatment and find cases where increasing the checkpoint inhibitor drugs decreases the efficacy of treatment. 

In this work, we explore a minimal oncolytic virus model to capture the essential features of virus infectivity in solid tumours and study how the dynamics of infectivity impact tumour outcome. Using a system of PDEs accounting for the behaviour of tumour cells and virus particles, we introduce an infectivity (or internalisation) parameter that regulates the ability of viral vectors to penetrate and infect tumour cells. This parameter is capable, to some extent, of accounting for the general penetrability in the tumour environment due to the extracellular matrix or variations in tumour cell receptors for virus internalisation such as lack of CAR expression. We do not explicitly account for virus clearance due to immune cells; however, we study changes in the behaviour of our system for variations in virus clearance rate. To understand how infectivity drives the behaviour in the spatiotemporal model, we also study the dynamics of the temporal model without diffusion by using bifurcation theory.\\

The work is organised as follows. In Sections~\ref{s1} and \ref{s2}, we introduce the PDE model and discuss results for tumours infected by viral agents with different infection rates. In Section~\ref{s3}, we derive and analyse a local ODE model for well-mixed populations that stems from the full PDE model. These equations underpin the behaviour of the full PDE model, and allow us to investigate the effect of changes in different parameters with respect to therapy outcomes, by means of bifurcation theory. In Section~\ref{s3b}, we compare the conclusions from the bifurcation analysis to the full PDE model.  In Section~\ref{s4}, we discuss the aspects that influence the quality of virotherapy as suggested by our findings. The paper is concluded in Section~\ref{s5}. Finally, \ref{appendix:parameterestimates} illustrates model calibration with respect to existing laboratory data, and \ref{appendix:parameterdependentequilibria} shows some analytical results.

\section{The model}\label{s1}
Let $ U(r,t), V(r,t)$ and $I(r,t)$ denote the population densities of uninfected tumour cells, free virus and infected tumour cells at distance $r$ from the centre of the tumour at time $t$. Assuming that diffusion is spherically symmetric in three dimensions, the equations that govern the interactions among populations are given by

\begin{align}
\frac{\partial U}{\partial t} & = \frac{D_u}{r^2} \frac{\partial}{\partial r} \left(r^2 \frac{\partial U}{\partial r} \right) + r_U U \left( 1- \frac{U+I}{k} \right) - \beta UV,\label{eq:TumourPDE} \\
\frac{\partial V}{\partial t} & = \frac{D_v}{r^2} \frac{\partial}{\partial r} \left(r^2 \frac{\partial V}{\partial r} \right) - \delta_v V - \beta (U + I) V + \alpha \delta_I I, \label{eq:VirusPDE}\\
\frac{\partial I}{\partial t} & = \frac{D_u}{r^2} \frac{\partial}{\partial r} \left(r^2 \frac{\partial I}{\partial r} \right) + \beta UV -  \delta_I I.\label{eq:InfectedPDE} 
\end{align}

\noindent In Eq.~(\ref{eq:TumourPDE}), diffusion and growth are modelled using Fisher's equation~\cite{fisher1937,murray2002} with growth being inversely proportional to tumour size. Here, $D_u$ is the diffusion rate of tumour cells, which is assumed constant, $r_U$ is the tumour's rate of growth and $k$ represents its carrying capacity. No vascular system is assumed to be present. The growth is hypothesised to be logistic and accounts for both uninfected and infected populations. The last term accounts for the infection of tumour cells by virus, which is modelled by mass action with infection rate $\beta$. 

Eq.~(\ref{eq:VirusPDE}) describes the dynamics of the virus, and the first term accounts for diffusion of viral particles with diffusivity $D_v$. The second term accounts for virus death at rate $ \delta_v $, and the third term describes internalisation of the virus, which, for this model, can occur in both uninfected and infected cells. Note that this term allows for loss of viruses to already infected cells through internalisation before lysis. For simplicity, the rate at which this takes place is also given by $\beta$, which is the same as for uninfected cells, $U(r,t)$. The final term describes the generation of new viral particles inside infected cells at the time of burst.  The parameter $\alpha$ is the mean number of new particles produced from each infected cell, or the viral burst size, and $\delta_I I$ is the death rate of infected cells. 

\begin{figure}[]
  \centering
  \begin{subfigure}[t]{0.03\textwidth}
    \textbf{A}
  \end{subfigure}
  \begin{subfigure}[t]{0.55\textwidth}
    \includegraphics[width=\linewidth, valign=t]{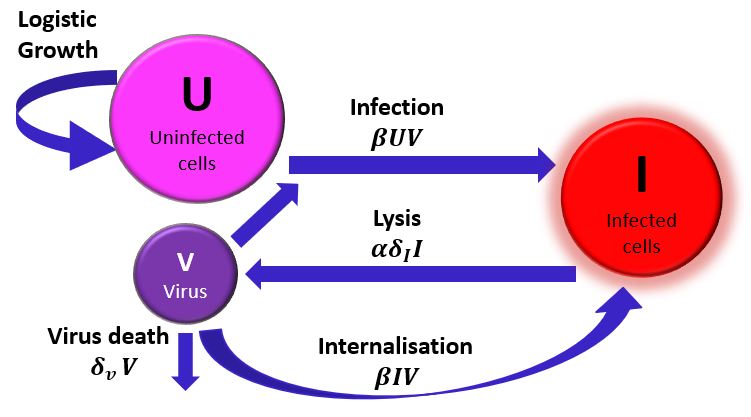}
  \end{subfigure}\hfill
  \begin{subfigure}[t]{0.03\textwidth}
    \textbf{B}
  \end{subfigure}
  \begin{subfigure}[t]{0.35\textwidth}
    \includegraphics[width=\linewidth, valign=t]{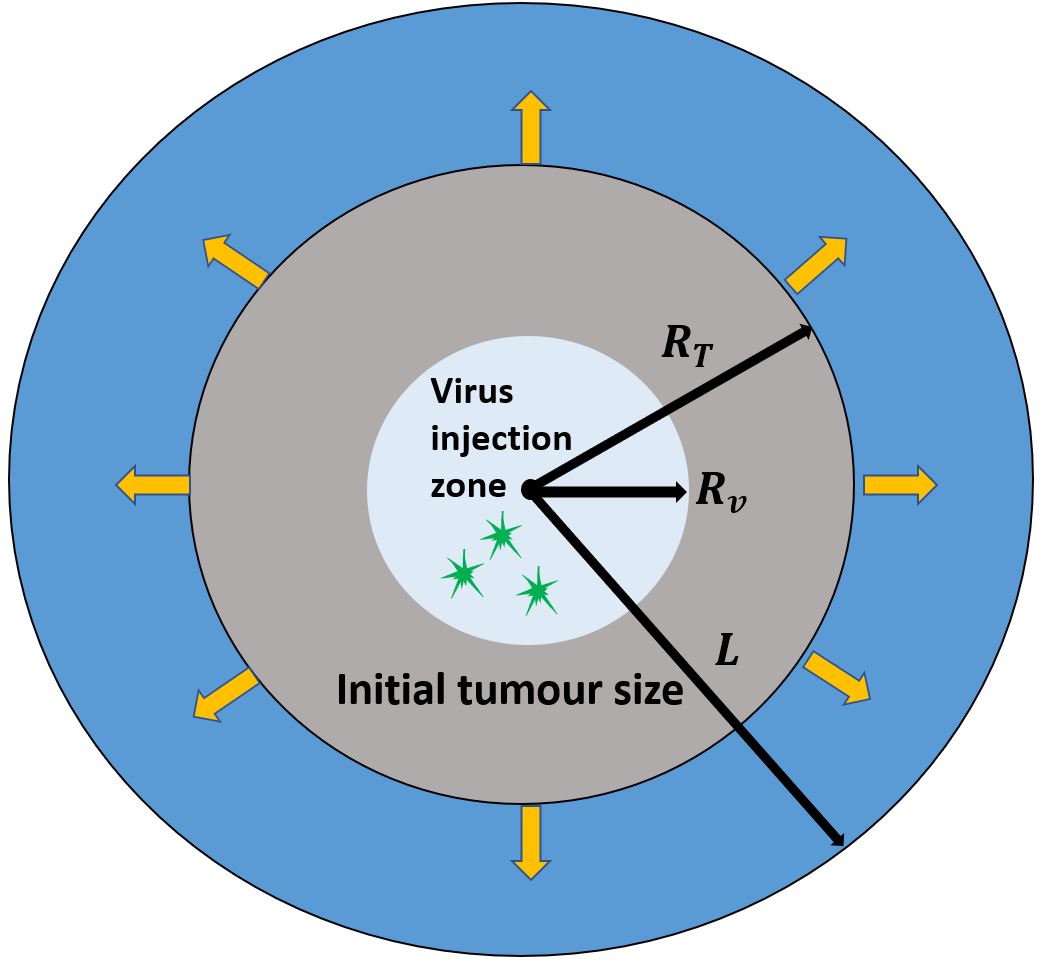}
  \end{subfigure}
  \caption[]{(A) Interactions among model populations. The uninfected tumour grows logistically at rate $r_U$. The viral population internalises and infects the uninfected tumour cells at rate $\beta UV$, and viral particles are lost to uninfected and infected cells at rates $\beta UV$ and $\beta IV$. The infected cells burst at rate $ \delta_I I$ and release $\alpha$ new virus particles per cell. Finally, the virus population decays at rate $\delta_v V$. (B) Initial conditions for the system (depicted in two dimensions only) with a tumour of radius $R_T$. The virus is injected in the centre of the tumour with radius of injection $R_V$, where $L$ is the maximum tumour radius. At the start, density is assumed equal to carrying capacity and the volume of tumour is set at $70$ mm$^3$.}
  \label{fig:flowchart}
\end{figure}

The infected population in Eq.~(\ref{eq:InfectedPDE}) is assumed to be diffusing at the same rate as uninfected cells $D_u$, whilst the remaining terms account for infection and death through lysis. All parameters in the model are assumed to be positive. Fig.~\ref{fig:flowchart}A  presents the interactions described by the system.

Initially, the spherical tumour has radius $R_T$ and the viral load is injected directly into its centre. The injection region is also spherical with radius $R_V$ (see Fig.~\ref{fig:flowchart}B), giving the following initial conditions
\begin{align} \label{Initial conds Inhomo}
& U(r,0) =
  \begin{cases}
    U_0      & \quad \text{if } 0 \leq r \leq R_T,\\
    0  & \quad \text{otherwise}, \\
  \end{cases}\\
& V(r,0) =
  \begin{cases}
    V_0       & \quad \text{if } 0 \leq r \leq R_V,\\ 
    0  & \quad \text{otherwise}, \\
  \end{cases}\\
& I(r,0) = 0,
\end{align}
where $U_0$ and $V_0$ are the initial densities of uninfected cells and virus. No-flux boundary conditions are assumed at both boundaries $r=0$ and $r=L$, where $L$ is the radius corresponding to the largest admissible tumour radius. Direct injection of virus into the tumour region has so far been a common practice in experimental settings. This is due to the antiviral immunity in the blood which will clear the virus population before they reach the tumour site~\cite{marchini2016overcoming, kim2015stem}. Experimental results indicate that the efficacy of intravenous delivery of oncolytic viruses is still limited and an area of focus in oncolytic virotherapy research~\cite{roy2013cell, seymour2016oncolytic}.

We calibrate and compare the dynamics from our model to those performed by Kim {\em et al.} (refer to Figure A in~\cite{kim2006relaxin}). This series of experiments investigates the role of the extracellular matrix in inhibiting viral spread. The authors track tumour volume in time between a standard adenovirus (Ad-$\Delta$E1b) and a relaxin-expressing adenovirus (Ad-$\Delta$E1b-RLX). The experiments are performed in vivo with nude mice using five different cell lines. The advantage in using these experimental results is that the primary function of relaxin is to degrade collagen and inhibit its synthesis so that the correlation between changes in tumour environment to treatment outcome are more pronounced.  

Parameter values for Eqs.~(\ref{eq:TumourPDE})--(\ref{eq:InfectedPDE}) and initial conditions are listed in Table~\ref{tab:Parameters}. The procedure used for calibrating the model is detailed in~\ref{appendix:parameterestimates}. The infectivity parameter $\beta$ is considered a free, phenomenological parameter: changes in its value describe cases where, for example, collagen has been degraded allowing easier access to tumour cell receptors or genetic enhancements to viruses to increase internalisation. The result is, on average, a variation in infiltration and internalisation patterns but, as we will see shortly, arbitrarily increasing $\beta$ does not always correspond to better outcomes or a faster tumour erosion overall.

\begin{table}
\centering
  \begin{tabular}{|c|c|c|c|}
    \hline
    \multicolumn{4}{|c|}{Parameter Estimates}\\
    \hline
    Symbol & Definition & Value & Reference \\
    \hline
    $k$ & Maximum tumour density & $10^6$ cells/mm$^3$ & \cite{lodish2008molecular}\\
    $r_U$ & Tumour growth rate & 0.3/day & \cite{cowley2014parallel}\\
    $D_v$ & Virus diffusion coefficient & 0.24 mm$^2$/day & \cite{paiva2009multiscale}\\
    $D_u$ & Tumour diffusion coefficient & 0.006 mm$^2$/day & \cite{kim2006relaxin}\\
    $\delta_v$ & Virus death rate & 4/day & \cite{mok2009mathematical}\\
    $\delta_I$ & Infected cell death rate & 1/day & \cite{ganly2000productive}\\
    $\beta$ & Internalisation rate constant & $\approx 1.5$ x $10^{-9}$ mm$^3$/ (viruses $\times$ day), variable & \cite{Fried2006}, estimated \\
    $\alpha$ & Viral burst size & 3500 viruses/cell & \cite{chen2001cv706}\\
    $U_0$ & Initial density of uninfected cells & $10^6$ cells/mm$^3$& estimated\\
    $V_0$ & Initial density of virus particles & 1.9 x 10$^{10}$ viruses/mm$^3$ & \cite{chen2001cv706}\\
    $R_T$ & Initial tumour radius & 2.6 mm & \cite{kim2006relaxin}\\
    $R_V$ & Initial radius of virus inoculation & 0.5 mm & estimated \\
    $L$ & Maximal tumour radius & 10 mm & estimated \\
    \hline
  \end{tabular}
  \caption{Parameters and initial conditions used in the model.}\label{tab:Parameters}
\end{table}

\begin{figure}[]
  \begin{subfigure}[t]{\textwidth}
    \textbf{A}
  \end{subfigure}
  \begin{subfigure}[t]{\textwidth}
  \centering
    \includegraphics[width=0.75\linewidth, valign=t]{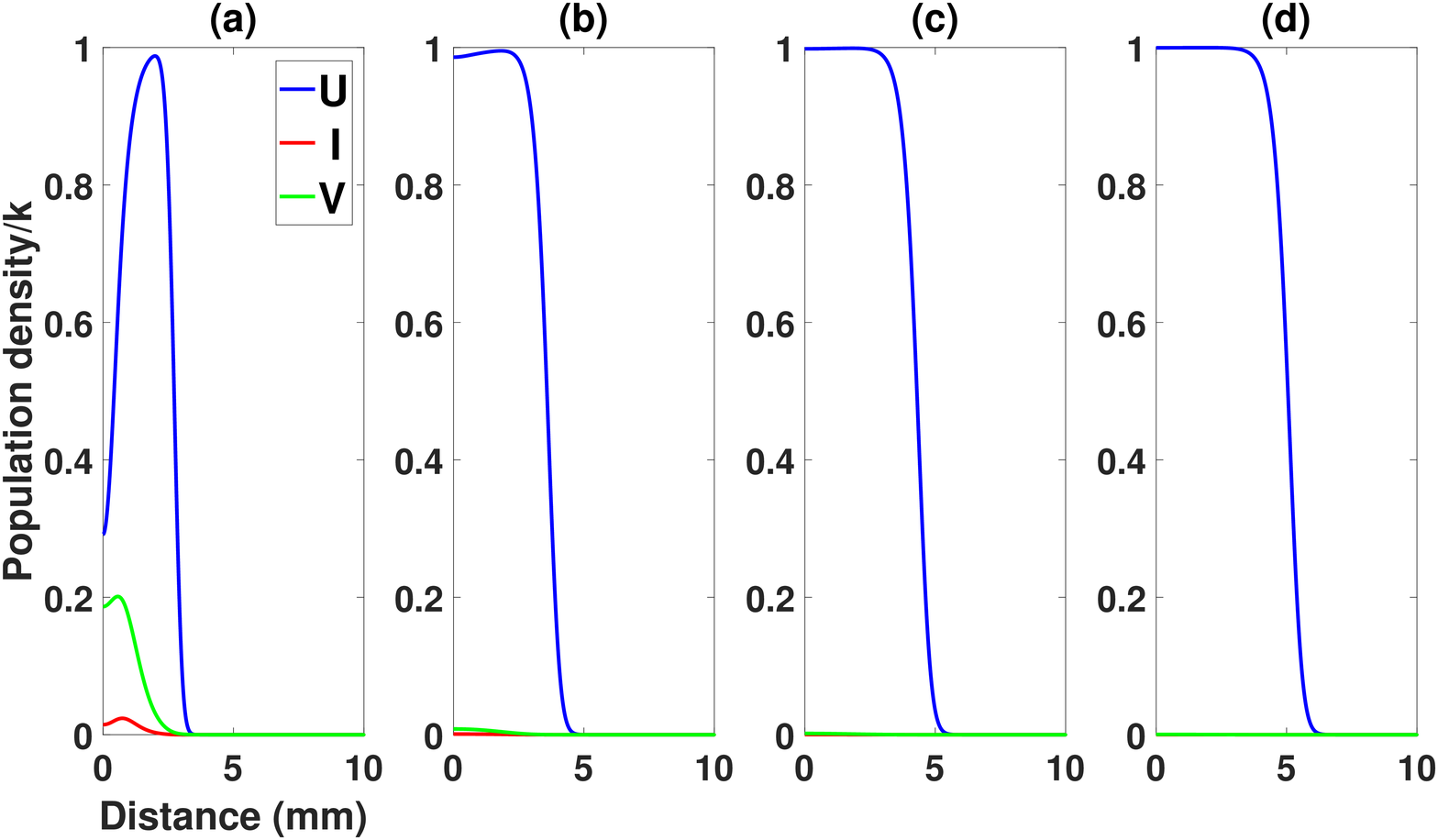}
  \end{subfigure}\hfill
  \begin{subfigure}[t]{\textwidth}
    \textbf{B}
  \end{subfigure}
  \begin{subfigure}[t]{\textwidth}
  \centering
    \includegraphics[width=0.75\linewidth, valign=t]{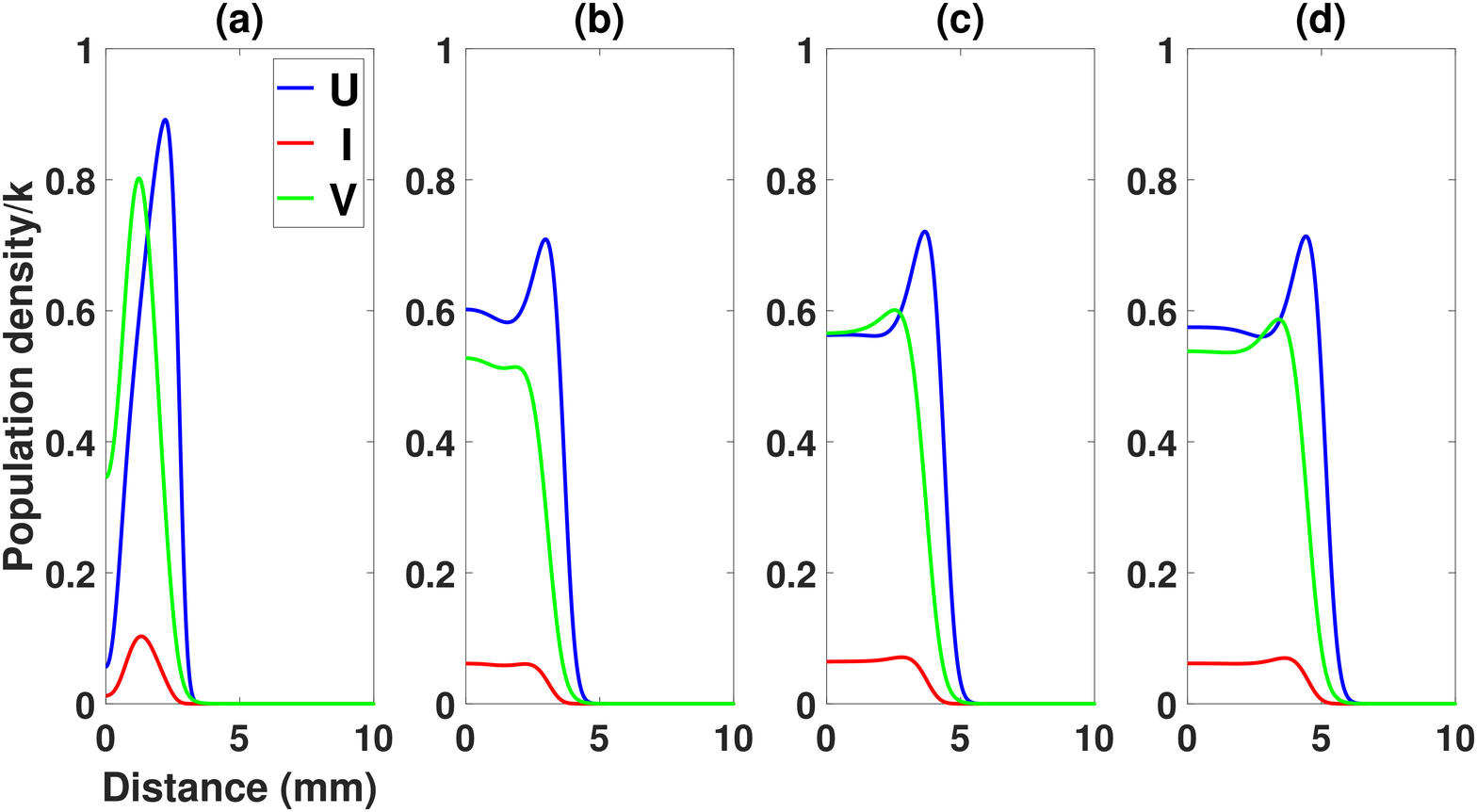}
  \end{subfigure}
  \begin{subfigure}[t]{\textwidth}
    \textbf{C}
  \end{subfigure}
  \begin{subfigure}[t]{\textwidth}
  \centering
    \includegraphics[width=0.75\linewidth, valign=t]{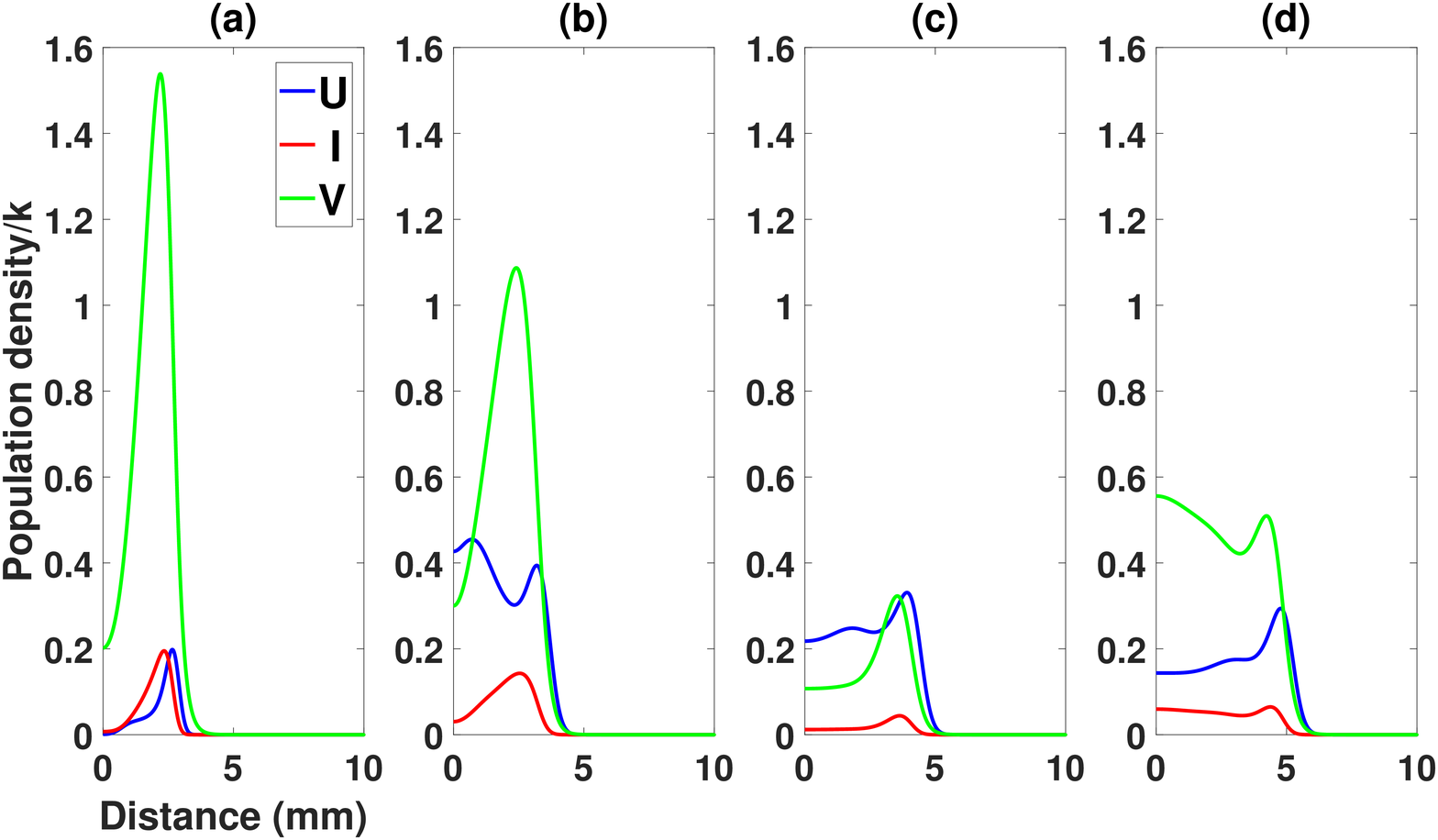}
  \end{subfigure}
  \caption[]{The evolution of population densities following viral injection for different values of $\beta$ at (a) $5$, (b) $20$, (c) $30$ and (d) $40$ days. Case (A) $\beta = 0.001$: after an initial response to viral treatment the virus population decays. Consequently, the tumour population recovers and reaches maximum density at $40$ days. Case (B) $\beta = 0.002$: there is a reduction in tumour density and waves appear stable at $40$ days, after minimal transient oscillations. Case (C) shows the largest $\beta = 0.005$: the tumour is considerably smaller after treatment, but oscillations in populations are larger and no stable waves are present at the end of the experiment. For readability, values of $V$ have been rescaled by $100$.}
  \label{fig:PDEdyn}
\end{figure}

\section{Results for localised tumours, for different infectivity}\label{s2}

Let us consider first the effect of oncolytic viral treatment in a dense tumour with low infectivity, monitoring the evolution of the system over $40$ days, which is the typical duration of the experiments from Kim {\em et al.} \cite{kim2006relaxin}. We assume a virus injection load of $V_0 = 1.9 \times 10^{10}$ viruses/mm$^3$, concentrated in a radial distance of $0.5$ mm around the centre of the tumour. As mentioned, variations in $\beta$ capture a relevant fraction of the changes in the tumour's environment or enhancements in virus binding with the overall goal of investigating changes in tumour outcome due to such modifications.

For reasons of simplicity, a scaled $\hat{\beta} = k\beta$ is considered in the following analysis, and we will drop the hat from now on. Starting at a low  $\beta = 0.001$, the treatment appears to have a minimal effect on the tumour population as shown in Fig.~\ref{fig:PDEdyn}A. Initially, the virus infects cells in the injection zone, reducing tumour density for uninfected cells (in blue) and increasing the numbers of infected cells (in red). After $5$ days, a wave of viral particles (in green) progresses towards the periphery, but the viral population already shows significant depletion in the area closest to the centre of the tumour. The infection rate is too weak to cause a sustained infection and the uninfected cells are close to carrying capacity after $20$ days: only few infected cells still survive and concentrate mostly around the tumour's centre. After $30$ days, there is no free virus left and no memory of the previous injection is present anywhere in the tumour. The uninfected population has fully recovered and continues growing as a standard travelling wave at maximum (carrying capacity) density. 

For an increased value of $\beta = 0.002$, corresponding to moderate infectivity, the dynamics are very different. At $5$ days, as shown in Fig.~\ref{fig:PDEdyn}B, viral particles reach higher densities close to the centre than in the previous example with consequent lower and higher maxima for uninfected and infected cells, respectively. In subsequent days, the populations experience a phase of dampened oscillations, after which the waves stabilise. At $20$ days, the uninfected cells are still being infected close to the centre and eventually reach a steady state around $30$ days since injection, after which their densities do not change with time. The travelling front of the wave of viral particles is responsible for an increase in infected cell numbers around the periphery, resulting from tumour expansion controlled by its growth and diffusion rates. At the end of the experiment, i.e. at $40$ days, the waves are stable. 

Considering a high $\beta = 0.005$, the evolution is characterised by faster infection and larger oscillations of tumour populations as Fig.~\ref{fig:PDEdyn}C reveals. In particular, the dynamics in the first $5$ days are different than the previous cases: there is an almost complete eradication of the tumour in the centre (blue curve) and almost equally dense populations of uninfected and infected cells in its vicinity. As infection and diffusion progress, oscillations become larger and last longer. Noticeably, the density of uninfected cells is almost four times lower than for $\beta=0.002$ (Fig.~\ref{fig:PDEdyn}B) within a $3$ mm radius from the centre. As time progresses, there are different profiles for the density of uninfected cells throughout the experiment with maxima and minima that tend to level out in the central area of the tumour as the experiment evolves, although there is no sign of stabilisation of waves at $40$ days. Population wave crests appear to be all concentrated near the periphery of the tumour at the end of the experiment.

\begin{figure}[]
  \centering
  \begin{subfigure}[t]{0.03\textwidth}
    \textbf{A}
  \end{subfigure}
  \begin{subfigure}[t]{0.45\textwidth}
    \includegraphics[width=\linewidth, valign=t]{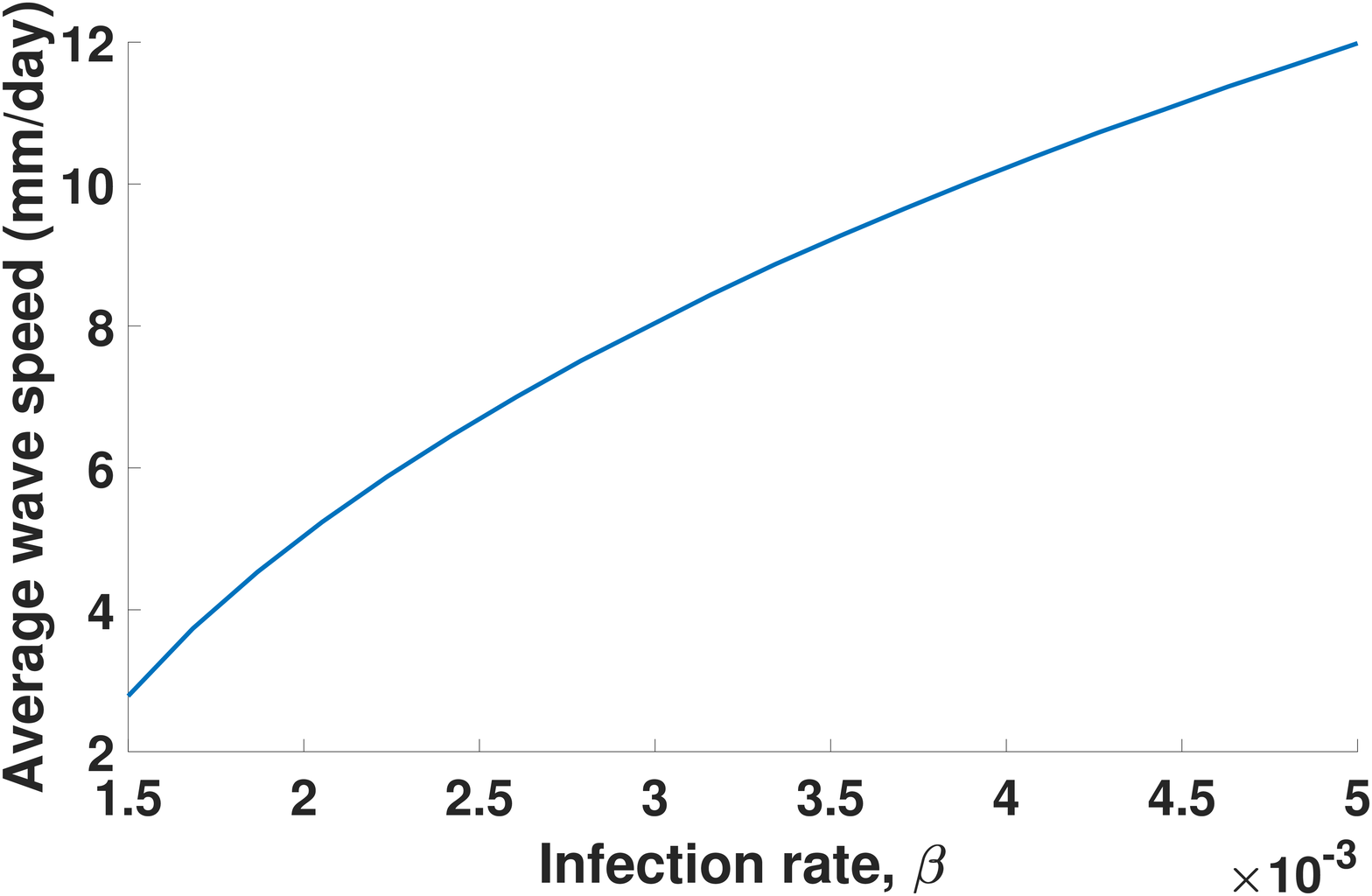}
  \end{subfigure}
  \begin{subfigure}[t]{0.03\textwidth}
    \textbf{B}
  \end{subfigure}
  \begin{subfigure}[t]{0.45\textwidth}
    \includegraphics[width=\linewidth, valign=t]{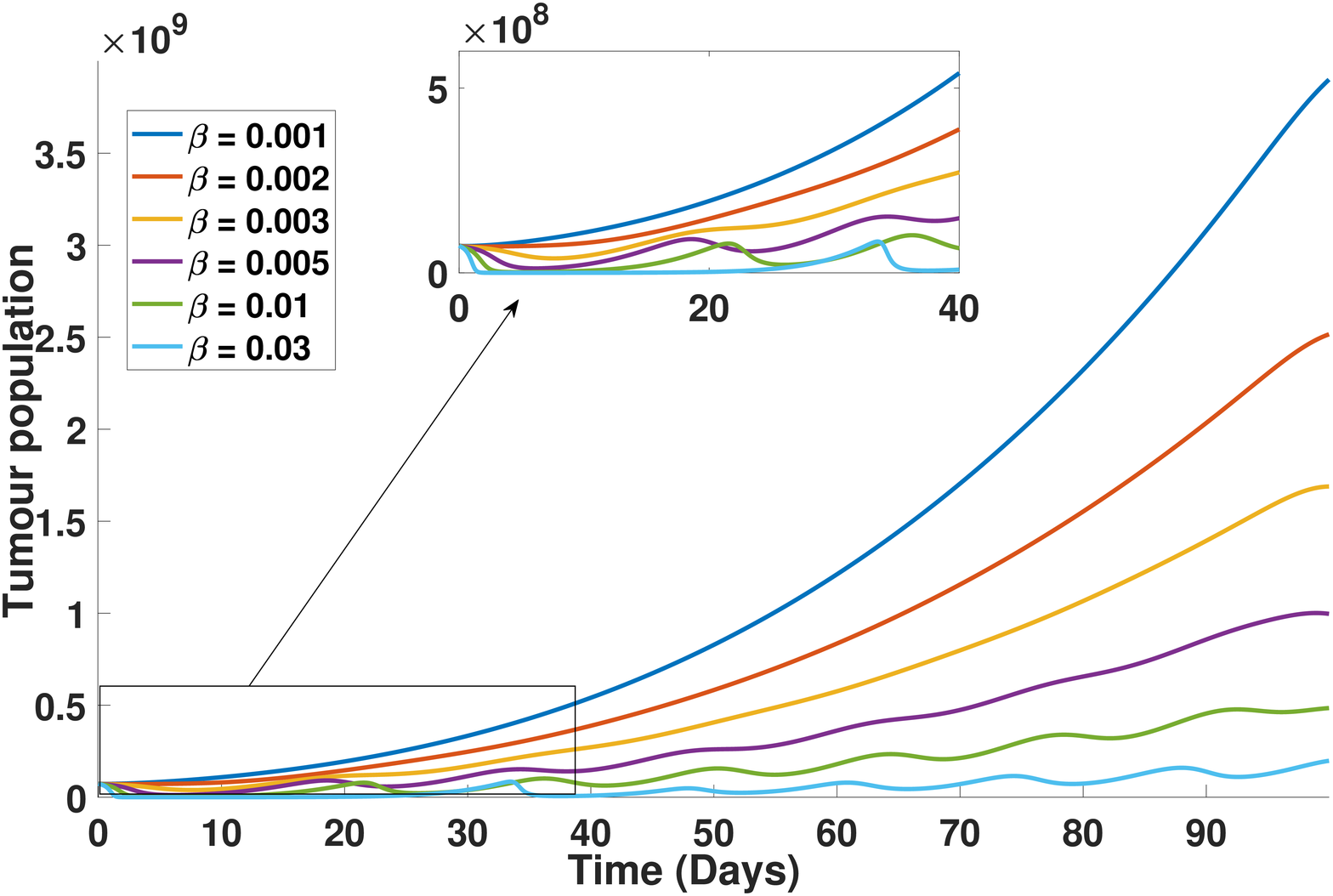}
  \end{subfigure}
  \begin{subfigure}[t]{0.03\textwidth}
    \textbf{C}
  \end{subfigure}
  \begin{subfigure}[t]{0.55\textwidth}
    \includegraphics[width=\linewidth, valign=t]{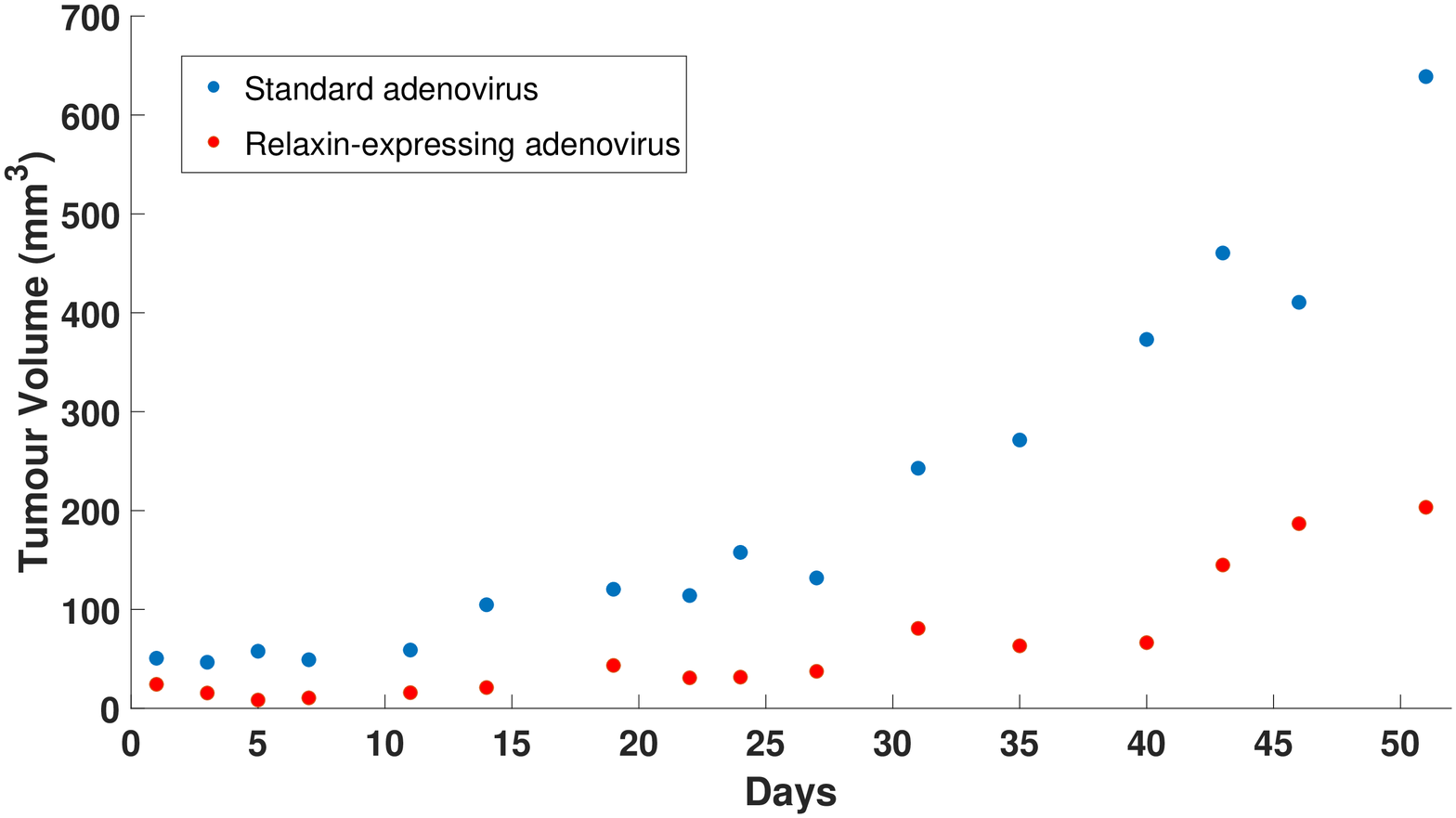}
  \end{subfigure}
  \caption[]{(A) Wave speed of virus propagation as a function of infection rate, $\beta$, using a linear interpolation scheme. The average wave speed of virus propagation through the tumour mass is depicted. As internalisation rate is increased, larger quantities of cells are infected which increases the wave speed. (B) Tumour population over $100$ days for different rates of $\beta$. Tumour decreases for larger $\beta$ values with noticeable oscillations appearing around $\beta \approx 0.003$.  The inset shows cell numbers in the first $40$ days of treatment. (C) Data from the Kim {\em et al.} experiments in~\cite{kim2006relaxin}. Growth of U343 cell lines in two mice. One mouse treated with a standard adenovirus (Ad-$\Delta$E1B) in blue and the other mouse treated with a relaxin-expressing virus (Ad-$\Delta$E1B-RLX) in red. Increased infectivity suppresses tumour size, and oscillatory behaviour is observable in the experimental results.}
  \label{fig:AvWaveSpeed}
\end{figure}

Within these values and as expected by intuition, a larger rate of internalisation and infection $\beta$ leads to better outcomes for growing tumours, but there are some important details to consider. Although Figure~\ref{fig:AvWaveSpeed}A shows an increase in the speed at which the virus propagates through the tumour without modifying the natural diffusion rate of the virus, a complete remission will require a very high infection rate, $\beta \approx 0.1$ (results not shown) and this may not be biologically possible. Note, since the system is continuous, we always have a resurgence of tumour cells; therefore, we assume tumour eradication when the total cell count is less than $1$. From the total cell count in time for increasing $\beta$, Figure~\ref{fig:AvWaveSpeed}B shows that reduction in tumour mass is always accompanied by oscillations in numbers of cancer cells, which occurs even at very high values for $\beta$. These oscillations start earlier in the treatment for higher internalisation rates and persist for large $t$. Also, the larger the $\beta$ the more delayed tumour growth on average is, because of viral killings.

It is important to note that the dynamics displayed by our model reproduces the dispersion and efficacy of viral treatment observed in available experimental results. Kim {\em et al.} measured the size of xenograft tumours in mice and observed that relaxin-expressing adenoviruses can delay tumour growth when compared with standard adenoviruses. Interestingly, their results are also oscillatory in nature as in Figure~\ref{fig:AvWaveSpeed}C, where the experimental results demonstrate that using the relaxin-expressing virus (in red) suppresses tumour growth as collagen is degraded in the tumour environment and viral infectivity is increased. Oscillations in tumour volume are also apparent in these results. Other experiments also depict oscillations where viruses are supported with collagen-degrading enzymes or antifibrotic agents, for example in Fig. $5$(A)-(B) in Guedan {\em et al.} \cite{guedan2010hyaluronidase} and Figure~$5$(A)-(B) in Diop-Frimpong {\em et al.} \cite{diop2011losartan}. For those cases, an increase in $\beta$ in our model seems to capture the salient features of those viruses' behaviour. The mathematical reason behind the birth of these oscillations will be clear in the next section.\\

\begin{figure}[]
  \centering
  \begin{subfigure}[t]{0.03\textwidth}
    \textbf{A}
  \end{subfigure}
  \begin{subfigure}[t]{0.45\textwidth}
    \includegraphics[width=\linewidth, valign=t]{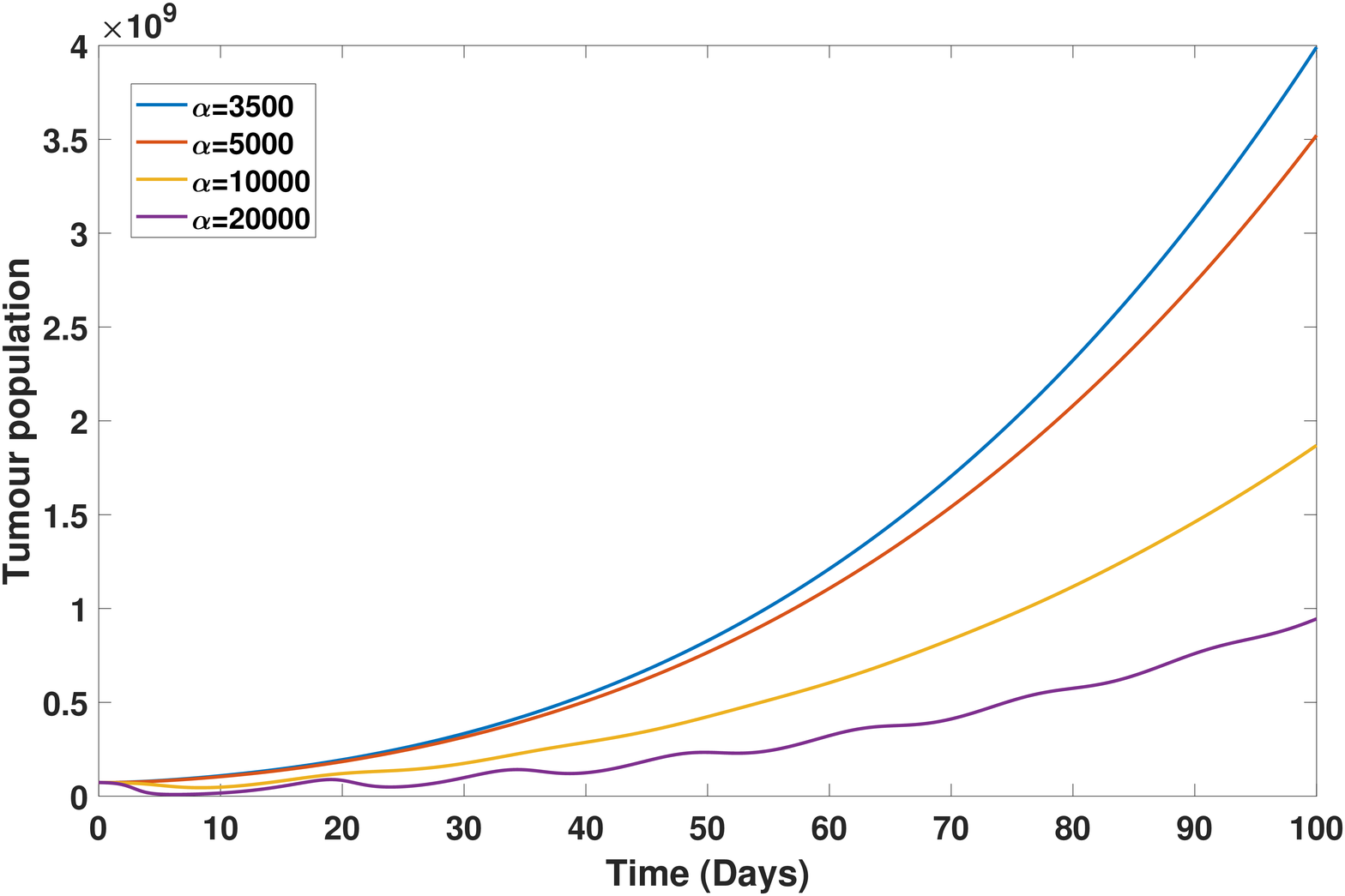}
  \end{subfigure}
  \begin{subfigure}[t]{0.03\textwidth}
    \textbf{B}
  \end{subfigure}
  \begin{subfigure}[t]{0.45\textwidth}
    \includegraphics[width=\linewidth, valign=t]{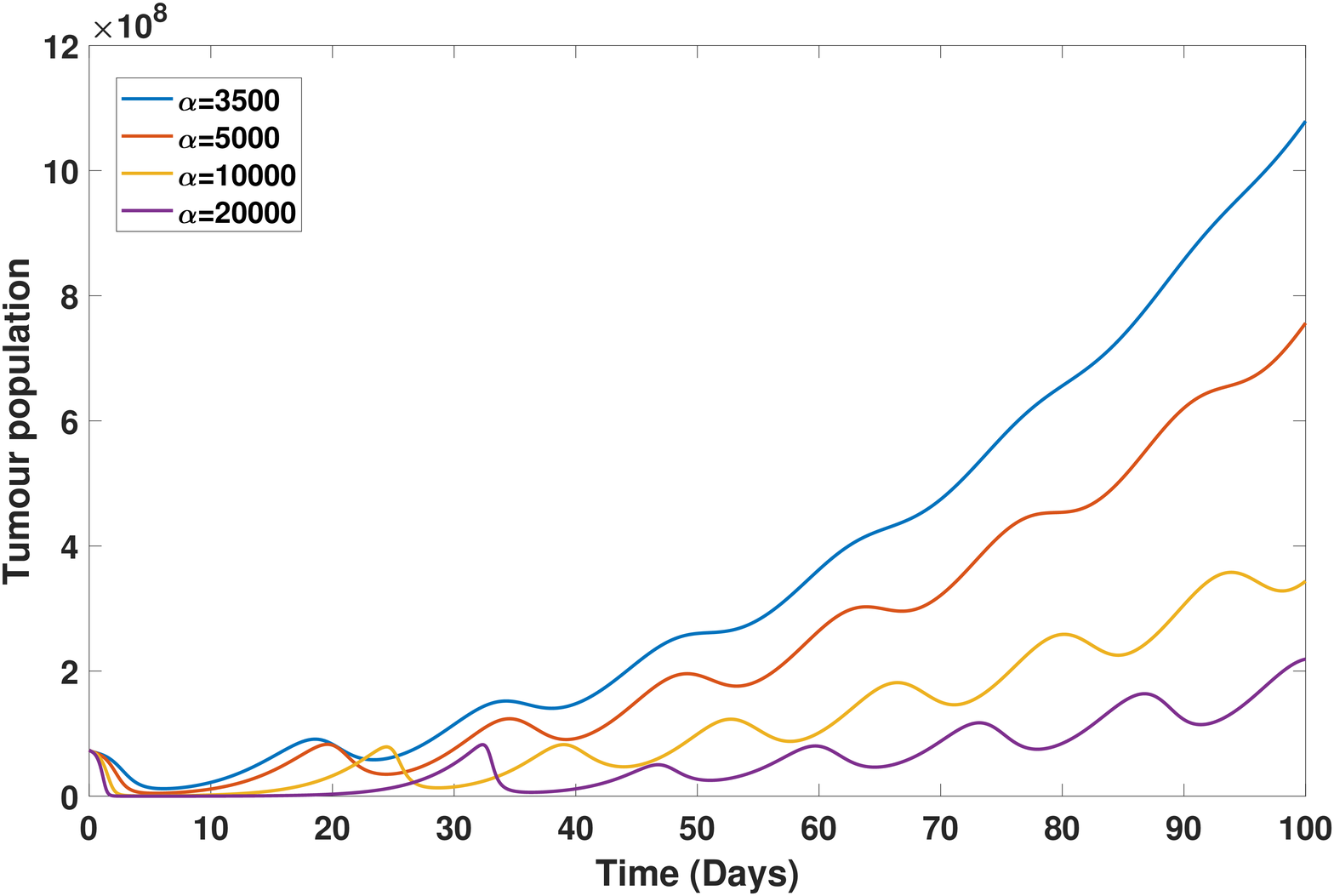}
  \end{subfigure}
  \caption[]{(A) Tumour population with low infectivity rate $\beta=0.001$ and varying proliferation rate $\alpha$. Only at very large $\alpha$ do we see a significant suppression of tumour volume. (B) At a high infectivity rate $\beta=0.005$, we see further reduction in tumour cell numbers as $\alpha$ is increased. Oscillations in tumour mass also increase in amplitude.}
  \label{fig:VaryingAlpha}
\end{figure}

We can also explore the changes in virus diffusion rate $D_v$ and proliferation rate $\alpha$. Increasing viral diffusion rate at low infectivity ($\beta=0.001$) or high infectivity ($\beta=0.005$) does not alter the dynamical behaviour of the system (results not shown). At low infection, this scenario could represent a tumour with reduced density but one which lacks the appropriate receptors for virus binding and internalisation so that the treatment remains ineffective. Increasing virus proliferation rate can reduce tumour size; however, this change needs to be relatively significant to impact cell numbers. At low $\beta$, as in Figure~\ref{fig:VaryingAlpha}A, increasing $\alpha$ from 3,500 (in blue) to 5,000 (in orange) does not significantly reduce tumour cell count. Even with $\alpha$ at 20,000 we cannot eliminate the tumour. We see the emergence of oscillations at high $\alpha$. At an increased $\beta$ of 0.005, increasing viral proliferation rate will result in a decrease in tumour load but we still do not have tumour clearance, as in Figure~\ref{fig:VaryingAlpha}B. At a viral burst size of 20,000, the minimum tumour cell count is 40,000 cells. Oscillations in tumour cell count also increase in amplitude as $\alpha$ increases. This result suggests that increasing virus reproductive capabilities, without enhancements to viral infectivity, may not alter the tumour outcome. In the next session, we will study the reduced ODE system to understand the appearance of these oscillations and to discover why tumour eradication is so difficult. We will also explore other viral enhancements through bifurcation analysis and compare the results from our temporal analysis to the full PDE model.

\section{The reduced ODE model, the origin of oscillations and virus design}\label{s3}

The theoretical impossibility for the PDE model to arrive at complete eradication can be traced back to mathematical properties of the defining equations. In particular, if homogeneous initial conditions in the injection zone are considered and we concentrate the analysis on that region, a well-mixed population of tumour cells and viral particles can be assumed. This allows an investigation into the stability and dynamical properties of a simplified system of ordinary differential equations (ODEs) once diffusion terms in the original PDE are set to zero. In this way, Eqs.~(\ref{eq:TumourPDE})--(\ref{eq:InfectedPDE}) reduce to the following nonlinear autonomous system, after the tumour populations are scaled by the carrying capacity $k$:
\begin{align}
\dot{U} & =  r_U U \left( 1- (U+I)\right) - \beta UV, \label{ode1} \\
\dot{V} & =  \alpha \delta_I I  - \delta_v V - \beta (U + I) V, \label{ode2} \\
\dot{I} & =  \beta UV - \delta_I I. \label{ode3}
\end{align}
This system shows typical features of existing models in the literature \cite{wodarz2001viruses,bajzer2008modeling,komarova2010ode,dingli2009dynamics,wodarz2003gene,jenner2018minimal} with some unique behaviours that stem from the depleting term for free viruses due to internalisation. In fact, the term $-\beta (U + I) V$ acts as a feedback on the viral load, which is reduced in proportion to the population of tumour cells, infected and uninfected, present in the system. As we will show shortly, this mechanism counteracts steady declines in tumour cells and prevents the system from arriving at a stable, tumour-free equilibrium.

Setting $ \dot{U} = \dot{V} = \dot{I} = 0$ for the equations above, two immediate solutions are given by $(U, V, I) = (0,0,0)$ and $(U, V, I) = (1,0,0)$. They are independent of parameters and represent complete eradication and a failed therapy attempt. The Jacobian matrix for Eqs.~(\ref{ode1})-(\ref{ode3}) is given by
\begin{equation}\label{JacobianODE}
JF = \left[
\begin{matrix}
r_U(1-I-2U)-V \beta & -U \beta & -r_U U \\
-V \beta & -(I + U) \beta - \delta_v & - V \beta + \alpha \delta_I \\
V \beta & U \beta & -\delta_I 
\end{matrix} \right].
\end{equation}
The equilibrium $(U, V, I) = (0,0,0)$ yields eigenvalues $\lambda_1 = r_U$ , $\lambda_2 = -\delta_I$ and $\lambda_3 = -\delta_V$. Because all parameters in the model are positive, $\lambda_1 > 0$ and the full-eradication solution is always unstable. This is why such a solution cannot be observed in the PDE version (\ref{eq:TumourPDE})-(\ref{eq:InfectedPDE}), independently of how fast the diffusive process or slow the tumour growth are. The equilibrium $(U, V, I) = (1,0,0)$ shows instead a dependence on system parameters, since its eigenvalues are given by
\begin{align}
\label{equilib2}
\begin{split}
\lambda_1 & = -r, \\
\lambda_{2,3} & = \frac{1}{2} \left( -(\beta + \delta_I + \delta_V) \pm \sqrt{( \beta + \delta_I + \delta_V)^2 - 4 \delta_I (\beta -\alpha \beta + \delta_V)}  \right).
\end{split}
\end{align} 
Terms $\lambda_1$ and $\lambda_3$ are always negative, and the stability of the system is determined by the sign of $\lambda_2$, which is zero for $\beta - \alpha\beta+\delta_V =0$ or 
\[
\beta^* = \frac{ \delta_v}{\alpha-1}.
\]
So, for $\beta > \beta^*$, the equilibrium associated to a failed treatment is unstable and, vice versa, it is stable for $\beta < \beta^*$ (as is the case, for example, for the values in Table~\ref{tab:Parameters} for the PDE system). As expected, the larger the decay rate $\delta_V$ of the virus particles, the larger the value of $\beta^*$ for which the therapy is ineffective. This means that, leaving all other model parameters unchanged, a shorter-lived virus requires a larger infectivity $\beta$ to have any reasonable impact on the tumour. Also, $\alpha$ is considered to be larger than one, since viral agents are expected to multiply inside infected cells.

Besides these immediate solutions, Eqs.~(\ref{ode1})-(\ref{ode3}) give rise to two further equilibria that are highly dependent on model parameters. It turns out that one of these solutions always occurs at negative values of the variable $U$ for biologically relevant parameter values and can thus be ignored. On the other hand, the other solution ${\bf U}_s=(U_s,V_s,I_s)$ can exist and be stable for values of interest, and its role for therapies with different values of $\beta$ will be clear shortly. For completion, we show analytical expressions for these equilibria in~\ref{appendix:parameterdependentequilibria}.

\begin{figure}[]
  \centering
  \begin{subfigure}[t]{0.02\textwidth}
    \textbf{A}
  \end{subfigure}
  \begin{subfigure}[t]{0.43\textwidth}
    \includegraphics[width=\linewidth, valign=t]{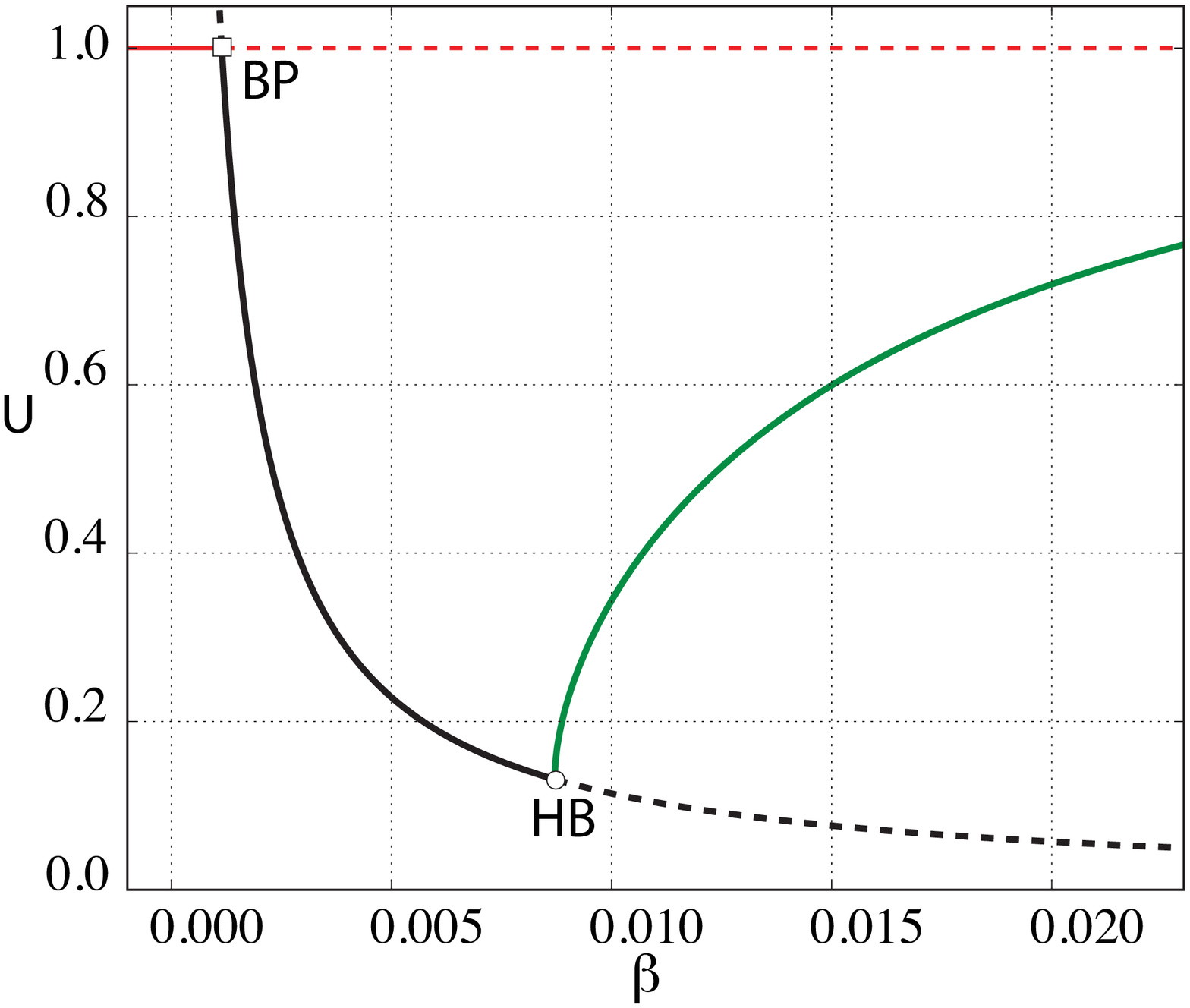}
  \end{subfigure}\hfill
  \begin{subfigure}[t]{0.02\textwidth}
    \textbf{B}
  \end{subfigure}
  \begin{subfigure}[t]{0.48\textwidth}
    \includegraphics[width=\linewidth, valign=t]{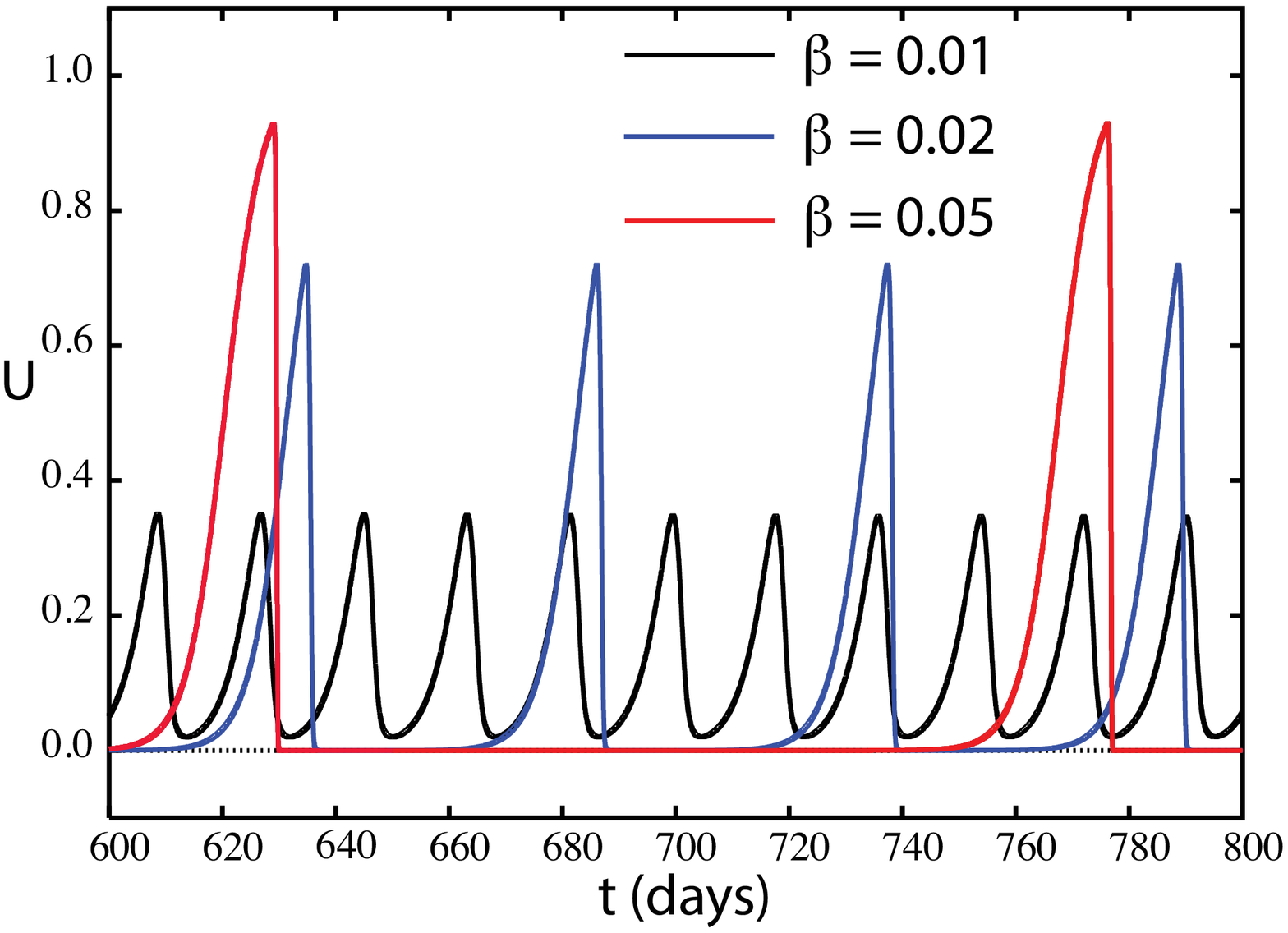}
  \end{subfigure}
  \caption[]{(A): One-parameter bifurcation plot for $\beta$ and the tumour population $U$ with other parameters as in Table~\ref{tab:Parameters}. Note that $U$ has been rescaled by the carrying capacity. The full-eradication solution at $U=0$ is always unstable (not shown). A branch point (BP) at $\beta_{\text BP} = 0.00114$ allows for a switch in stability for the ineffective treatment at $U=1$ with the emergence of a stable solution whose $U$ value decreases with $\beta$. A second bifurcation (HB) at $\beta_{\text HB} = 0.00871$ changes the stability of this solution and allows oscillations to emerge. This branch (in green) is depicted by using only the maxima of oscillations and is present for all values $\beta > \beta_{\text HB} $. Values of maxima increase towards $U=1$ as $\beta$ increases to large, biologically unrealistic values.  In panel (B), a sample of these orbits is shown for high values of $\beta$. As $\beta$ further increases, maxima and periods of oscillations increase, tending to the limiting case $U_{\text max} = 1$ for extremely large values of infectivity.}
  \label{fig:Bif_beta}
\end{figure}

A typical bifurcation plot of solutions in relation to infectivity $\beta$ is illustrated in Fig.~\ref{fig:Bif_beta}A for parameter values in Table~\ref{tab:Parameters}. Stable branches are indicated as continuous lines, whereas unstable ones are dashed. At very small values of $\beta$, i.e. $0<\beta<0.00114=\beta_{\text BP}$, tumour and virus characteristics do not allow for any effect on the virotherapy with the solution for uninfected cells at carrying capacity being stable.  At $\beta = \beta_{\text BP}$, a branch point (BP) occurs through which the original stable ineffective solution loses stability and allows for a second, stable branch to exist. This branch, whose value for the uninfected cells monotonically decreases for increasing values of $\beta$, corresponds to the aforementioned equilibrium solution ${\bf U}_s$, for which cells and virus remain at constant values. For $\beta < \beta_{\text BP}$, ${\bf U}_s$ is stable only at negative, non-biological values of $V$ and $I$ (not shown in the plot). As seen in another similar ODE model with different growth dynamics for cancer \cite{jenner2019gomp}, the presence of ${\bf U}_s$ for biologically meaningful values of internalisation shows that the therapy is still able to reduce the tumour, but a full eradication is not possible. We indicate this outcome as a co-existing equilibrium. As infectivity grows, the branch eventually undergoes a supercritical Hopf bifurcation (HB) at $\beta_{\text HB} = 0.00871$, where oscillations are born. These oscillations persist for all values of $\beta>\beta_{\text HB}$, with the limit cycle branch growing steadily to higher values of $U$. As infectivity increases, the maximum and minimum of the limit cycles tend to one and zero, respectively, as shown in Fig.~\ref{fig:Bif_beta}B. There is no intersection between the oscillation and the steady-state solution branches as well as no branching points or other types of bifurcations that can change their stability at finite values of $\beta$. The character of limit cycles also implies that, although a stable full-eradication solution is absent, for very large, biologically inadmissible values of $\beta$, oscillations show greater period and larger interval in between ``spike''-like orbits. 

A very similar scenario exists for $\alpha$, showing that the ability of virus to reproduce after infecting tumour cells plays a role in the success of virotherapy as well. For example, for a high level of infectivity $\beta=0.005$, we found that an increase in reproductive ratio has the ability to decrease the tumour density of uninfected cells. The shape of the bifurcation plot (not shown here) for $\alpha$ is the same as that for $\beta$ with a BP and HB respectively at $\alpha_{\text BP} = 800$ and $\alpha_{\text HB} = 6095$ with similar limit cycles that emerge and persist for all larger values of $\alpha$. The existence of a BP means that there is always a minimal amount of virus that has to be produced when cells lyse for the therapy to yield any benefit. In general, further bifurcation plots (not shown here) confirm that the lower the $\alpha$, the higher the $\beta$ required for viral loads to have any effect on tumours as expected.

\begin{figure}[]
  \centering
  \begin{subfigure}[t]{0.02\textwidth}
    \textbf{A}
  \end{subfigure}
  \begin{subfigure}[t]{0.47\textwidth}
    \includegraphics[width=\linewidth, valign=t]{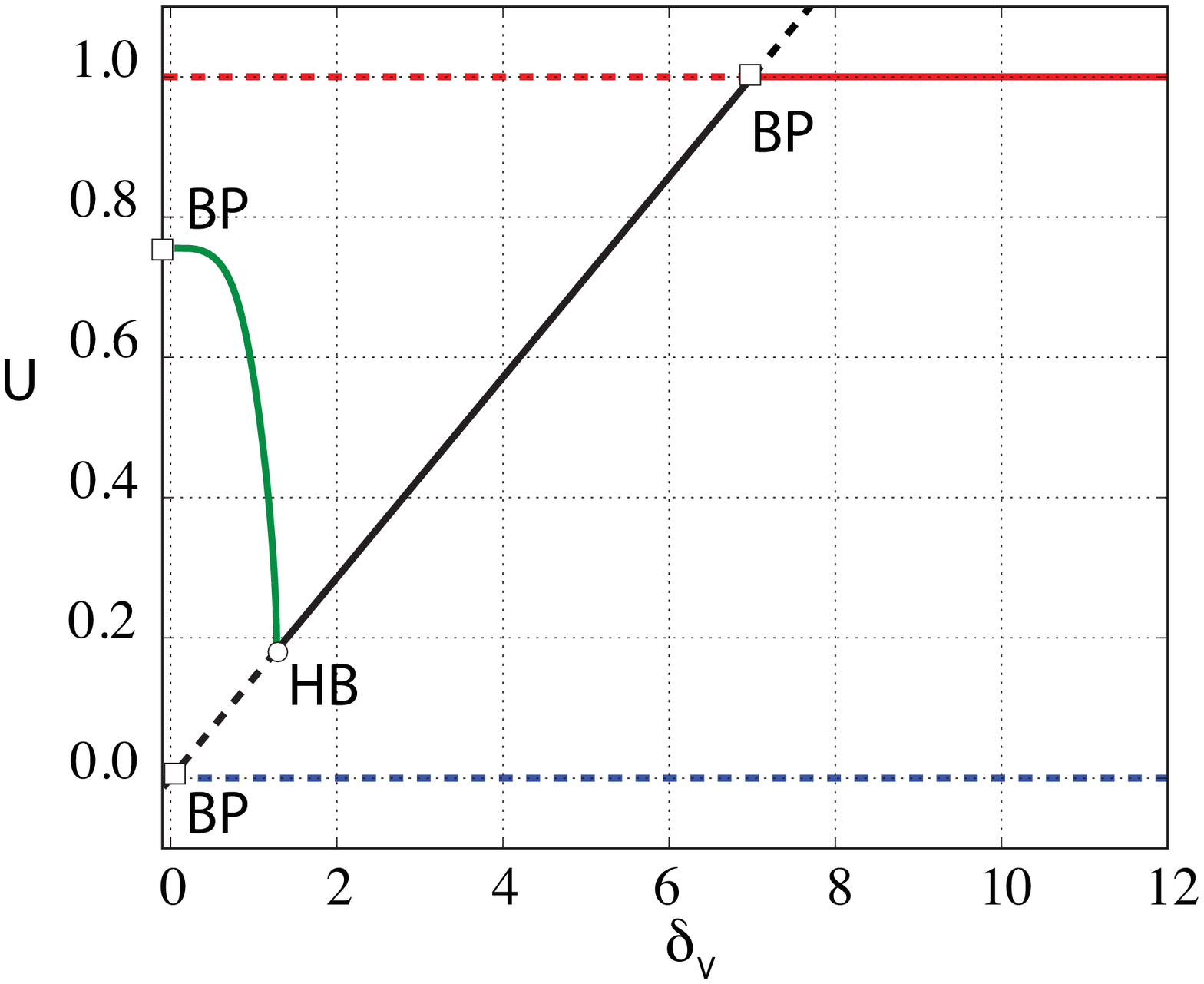}
  \end{subfigure}\hfill
  \begin{subfigure}[t]{0.02\textwidth}
    \textbf{B}
  \end{subfigure}
  \begin{subfigure}[t]{0.45\textwidth}
    \includegraphics[width=\linewidth, valign=t]{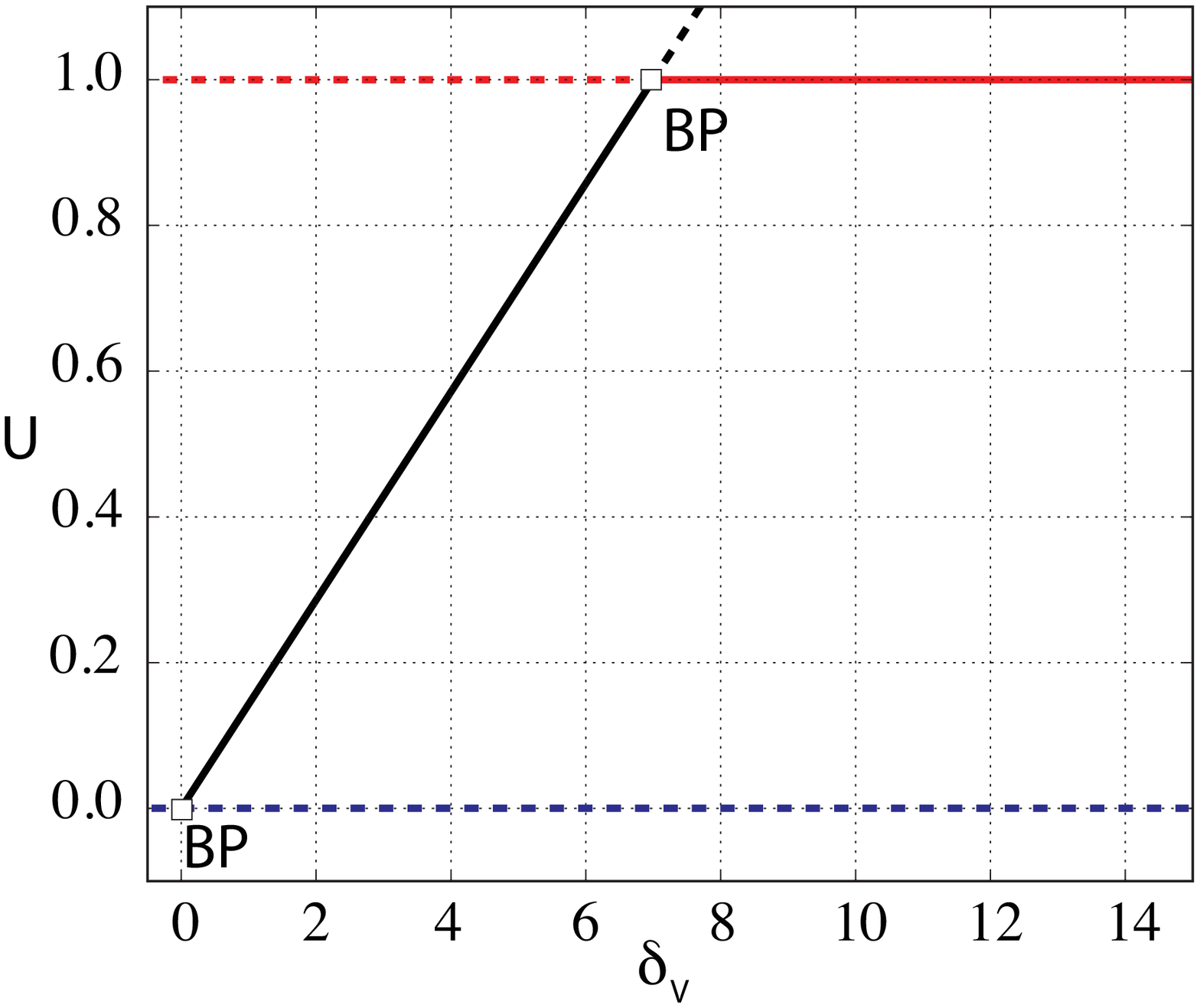}
  \end{subfigure}
  \caption[]{One-parameter continuations showing the role of the death rate constant for virus particles $\delta_V$ at an intermediate rate of infectivity $\beta=0.002$. (A) For $\delta_I = 1.2$, one Hopf point is present for $\delta_V = 1.276$ and limit cycles originating from it have an increasing maximum as $\delta_V$ is reduced. Oscillations terminate at a branch point (BP) at the limiting value $\delta_V = 0$. (B) For a high $\delta_I = 9$, no oscillations are possible and the co-existing equilibrium is present for all $\delta_V < \delta_{\text BP} = 6.998$. The value of $U$ for the co-existing equilibrium linearly decreases with $\delta_V$ and reduces to full eradication in the limiting, biologically meaningless case of an ``immortal'' virus, i.e. $\delta_V = 0$. }
  \label{fig:dr_b}
\end{figure}

Different and less intuitive is the role of death rate constants for viruses and infected cells. In Fig.~\ref{fig:dr_b}A, a bifurcation plot for $\beta=0.002$ and $\delta_I=1.2$ is shown with other parameters as in Table~\ref{tab:Parameters}. Interestingly, a HB emerges at low values of $\delta_V$ and stable oscillations (in green) grow their maxima steadily as $\delta_V$ further decreases.  This implies that, at moderate infectivity and infected cell death rates, extending the lifespan of the virus can result in oscillations rather than reduction of the tumour mass with periodic growths to densities that are higher in value than the co-existing equilibrium (in black). Nonetheless, low values of viral death rate reduce the value of tumour load $U$ at the co-exisiting equilibrium (compare it with Fig.~\ref{fig:Bif_beta}A). In the case of sufficiently large $\delta_I$, a complete-eradication solution exists for $\delta_V = 0$, with no presence of oscillations, as shown in Fig.~\ref{fig:dr_b}B. This is of course an idealised situation, where the virus does not decay and is able to survive (and attack) forever.

\begin{figure}[]
  \centering
  \begin{subfigure}[t]{0.02\textwidth}
    \textbf{A}
  \end{subfigure}
  \begin{subfigure}[t]{0.46\textwidth}
    \includegraphics[width=\linewidth, valign=t]{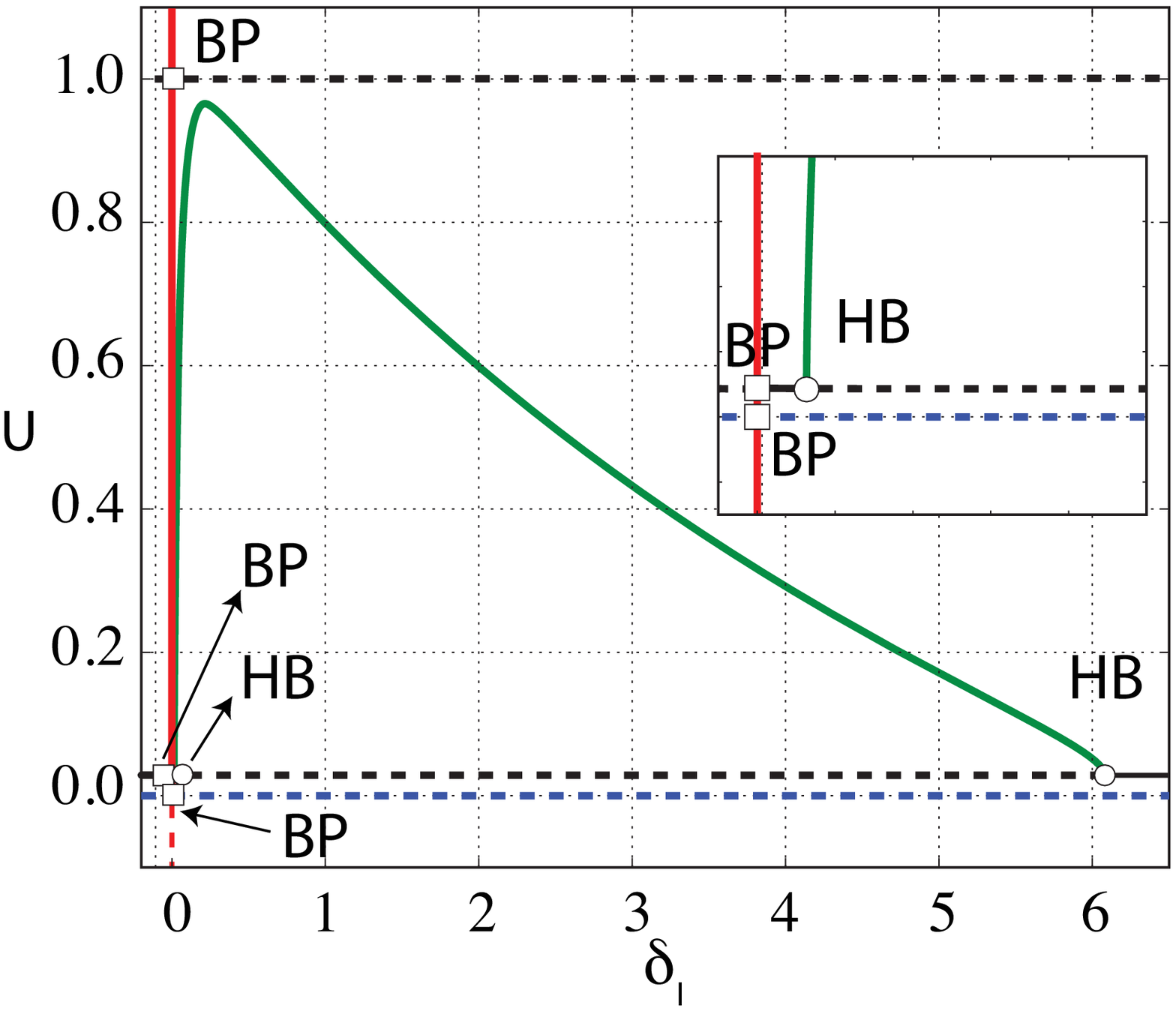}
  \end{subfigure}\hfill
  \begin{subfigure}[t]{0.02\textwidth}
    \textbf{B}
  \end{subfigure}
  \begin{subfigure}[t]{0.46\textwidth}
    \includegraphics[width=\linewidth, valign=t]{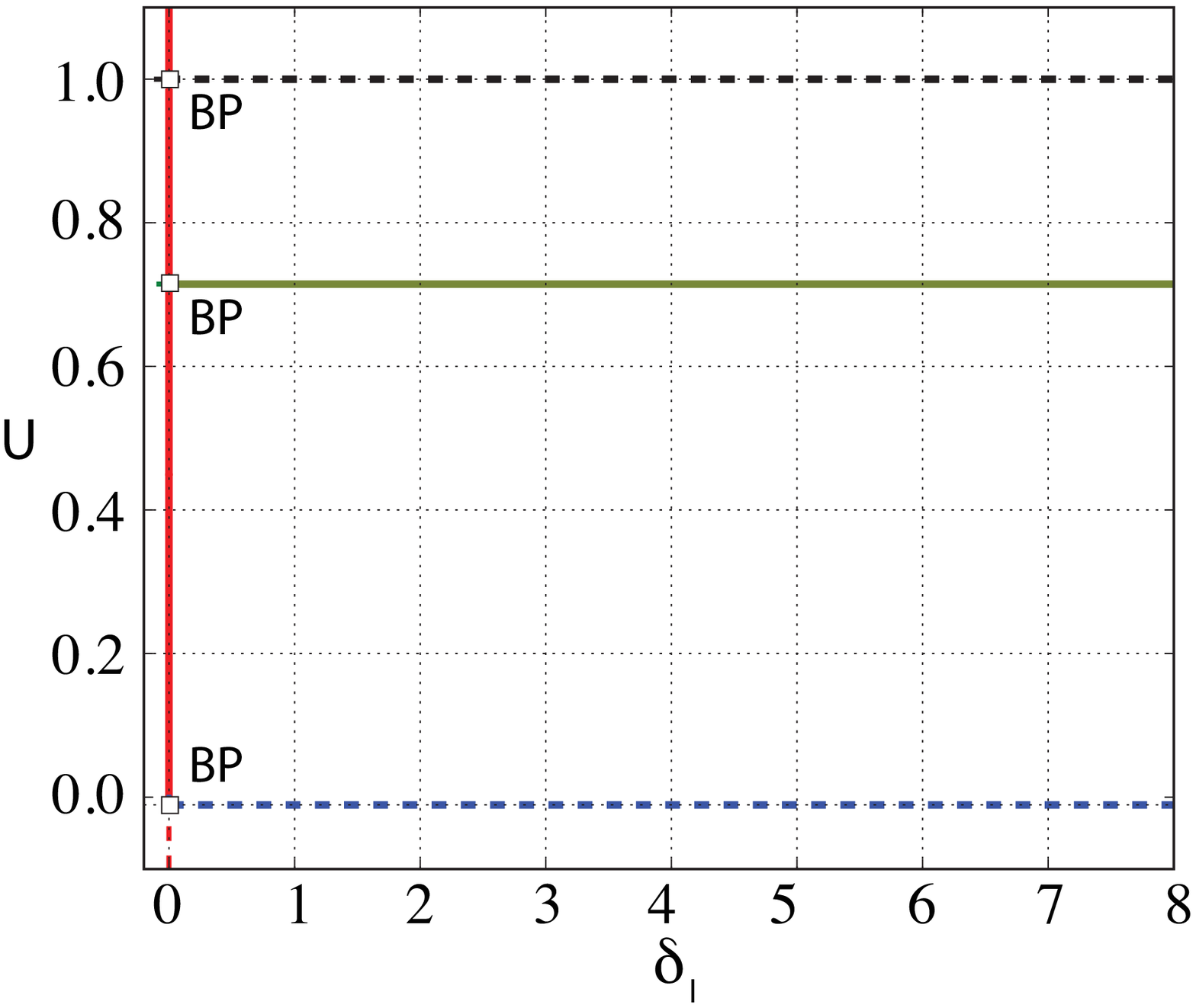}
  \end{subfigure}
  \caption[]{One-parameter continuations showing the role of the death rate constant for infected cells $\delta_I$, again at $\beta=0.002$. (A) For very low $\delta_V = 0.2$, two Hopf points connected by a branch of stable limit cycles (in green) are evident at $\delta_I = 0.0126$ and $\delta_I = 6.081$. The inset shows the bifurcation points occurring at values of $\delta_I$ close to zero. A further stable branch (in red) for $\delta_I = 0$ also exists. (B) For a higher $\delta_V = 5$, a co-existence equilibrium (in green) is present and is unaffected by changes in $\delta_I$. The same two branch points as in panel (A) for $\delta_V = 0$ corresponding to the birth of unstable solutions are present. A third branch point now exists at the intersection between the co-existence branch (in green) and the stable branch at $\delta_I = 0$.}
 \label{fig:dr_c}
\end{figure}

The effect of cell death rate $\delta_I$ on therapies is also worth considering. In Fig.~\ref{fig:dr_c}A, for a very low virus death rate $\delta_V = 0.2$, a faster or slower death rate for infected cells does not change the value of tumour load $U$ for the co-existing equilibrium (in black). Two Hopf points instead appear and, for biologically admissible values, a branch of stable limit cycles (in green) connects them. Note how there is a value for $\delta_I$ where the maximum of the oscillations in $U$ is highest and very close to carrying capacity. This is somewhat unexpected and counterintuitive, because it indicates that, depending on the other parameters, there is an intermediate death rate for infected cells $\delta_I$ that can create very high oscillations in tumour loads, notwithstanding a significant increase in virus lifespan (i.e. a very low $\delta_V = 0.2$). At larger and biologically sound values for $\delta_V$ (Fig.~\ref{fig:dr_c}B), a co-existing equilibrium (in green) appears and is unaffected by changes in $\delta_I$. Note also in both panels that, for the chosen parameter values, the equilibria corresponding to full eradication (in blue) and ineffectiveness of therapy (in black) are always unstable. A new, limiting and nonbiological equilibrium corresponding to the case of $\delta_I=0$ is present (in red), indicating that the number of initially uninfected cells $U$ remains constant because infected cells never die.

\begin{figure}
    \centering
    \includegraphics[width = 0.8\textwidth]{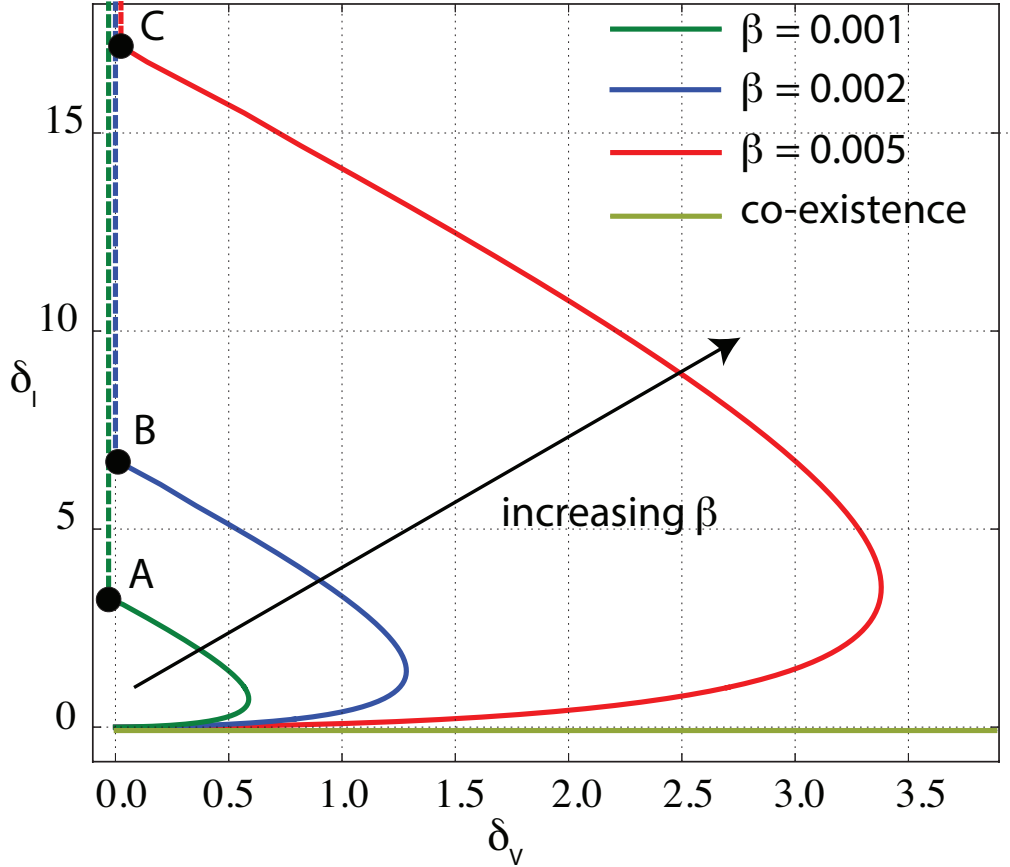}
    \caption{Branches of Hopf bifurcations in the $\delta_V$-$\delta_I$ parameter plane for different values of $\beta$. These curves  and the $\delta_I$-axis delimit the areas where oscillations can occur, with no limit cycles extending outside their boundaries. Limiting solutions at $\delta_V=0$ and $\delta_I = 0$ correspond to eradication and co-existence, respectively. Although co-existence (in olive) is possible for any value of $\delta_V > 0$ (i.e. it corresponds to the positive $\delta_V$-semi-axis), eradication only occurs for values $\delta_I$ that are larger than the ordinate of the intersection point for a given Hopf curve and the axis $\delta_V=0$ at a given $\beta$. These points are indicated with letters A, B and C respectively for $\beta$ equal to $0.001, 0.002$ and $0.005$. The dashed lines indicate where eradication occurs with colours referring to the given $\beta$. Note that, for readability, eradication lines that lie on the $\delta_I$-axis have been slightly shifted to the left so that they do not overlap.}
    \label{fig:dr_2b}
\end{figure}

The birth and extension of oscillations in tumour densities are driven by different parameters. In Fig.~\ref{fig:dr_2b}, branches of Hopf points are continued in the two parameters controlling virus and infected cell decay. Curves at different $\beta$ are presented, and the area delimited by these curves and the segment of the $\delta_I$-axis in between represent the part of parameter space where limit cycles are occurring. Larger infectivity increases the area of such oscillations, because branches of Hopf are more extended for larger $\beta$.  For values $k_1$ and $k_2$ such that lines $\delta_V = k_1$ or $\delta_I = k_2$ do not intersect the curves in Fig.~\ref{fig:dr_2b}, no Hopf points are encountered and no oscillations are present when the parameter $\delta_V$ or $\delta_I$ is varied, respectively. Note also that there exist two limiting cases of interest, also encountered when we have previously discussed the behaviour of one-parameter bifurcations. If $\delta_I = 0$, for any value of $\delta_V > 0$ a co-existing solution is expected, as shown in the inset of Fig.~\ref{fig:dr_c}A by the stable branch in between the BP and HB (in black), or as indicated in Fig.~\ref{fig:dr_c}B by the stable, unchanging vertical branch (in red). These occurrences are indicated in Fig.~\ref{fig:dr_2b} by a continuous olive line. If $\delta_V=0$, eradication is also possible as shown previously in Fig.~\ref{fig:dr_b}B, but this happens only if $\delta_I$ is sufficiently large. In fact, the value of $\delta_I$ has to be larger than the $\delta_I$-coordinate of the highest point of intersection between the Hopf branches and the axis $\delta_V=0$. This value indicates a threshold value for the decay rate of infected cells for which, in the limiting case of an immortal virus, the tumour can be completely destroyed.  In fact, analytically, when $\delta_V =0$, the Jacobian for Eqs.~(\ref{ode1})-(\ref{ode3}) (see ~\ref{appendix:parameterdependentequilibria}) associated to the equilibrium solution $(0,0,0)$ bears one negative and two zero eigenvalues, making a full eradication solution possible. This solution is, of course, not of biological interest but still indicates that, within the model, an ``invincible'' virus is not sufficient to clear the tumour if infected cells are not dying quickly enough. Given the strong analogy between the role of $\alpha$ and $\beta$ as previously discussed, we finally note that similar curves exist (not shown) also for increasing $\alpha$'s, with the area of limit cycles growing as $\alpha$ grows.

\section{Comparing the spatio-temporal PDE model against the bifurcation analysis of the ODE model}\label{s3b}
\begin{figure}[]
  \centering
  \begin{subfigure}[t]{0.02\textwidth}
    \textbf{A}
  \end{subfigure}
  \begin{subfigure}[t]{0.46\textwidth}
    \includegraphics[width=\linewidth, valign=t]{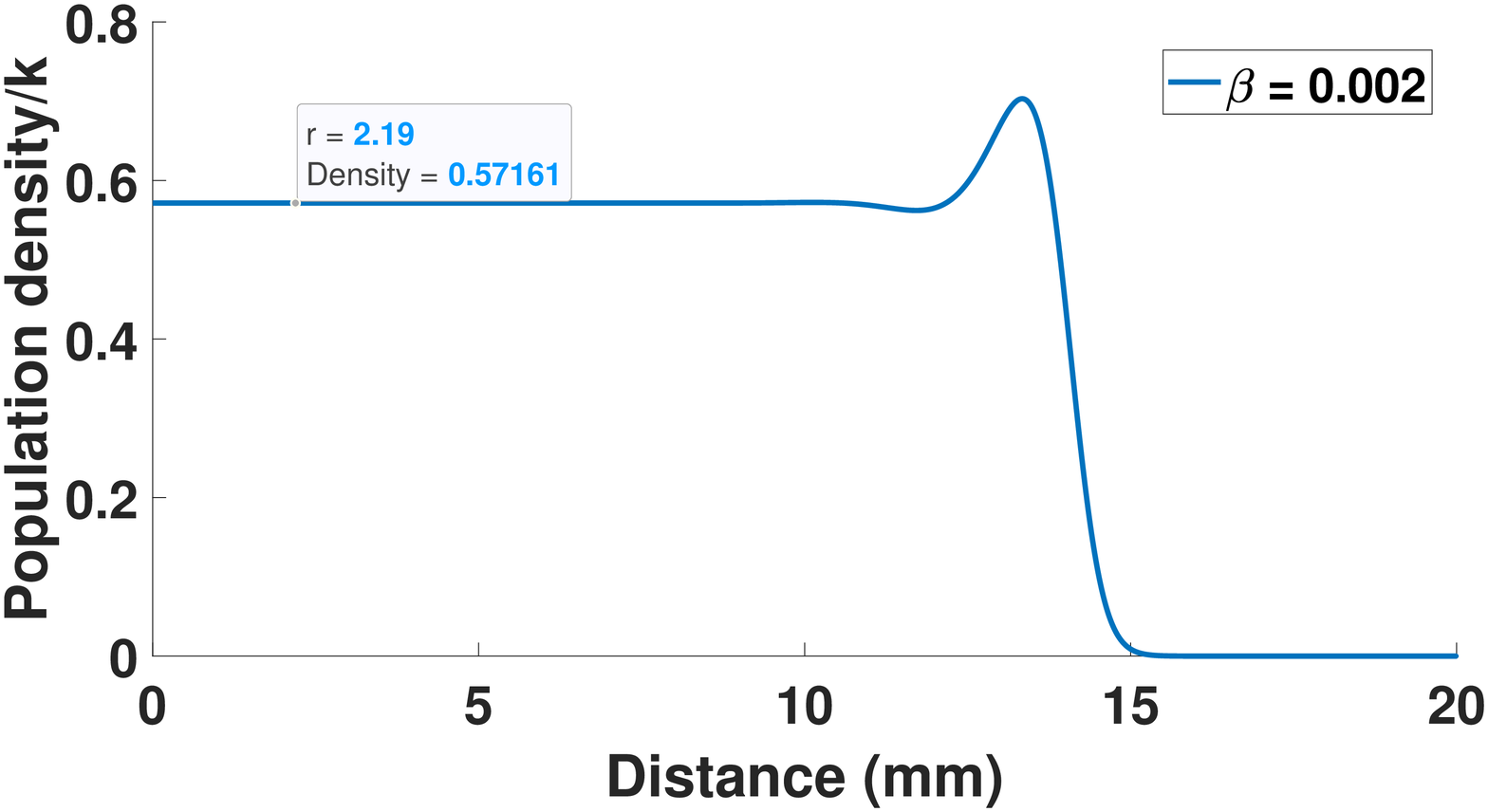}
  \end{subfigure}\hfill
  \begin{subfigure}[t]{0.02\textwidth}
    \textbf{B}
  \end{subfigure}
  \begin{subfigure}[t]{0.46\textwidth}
    \includegraphics[width=\linewidth, valign=t]{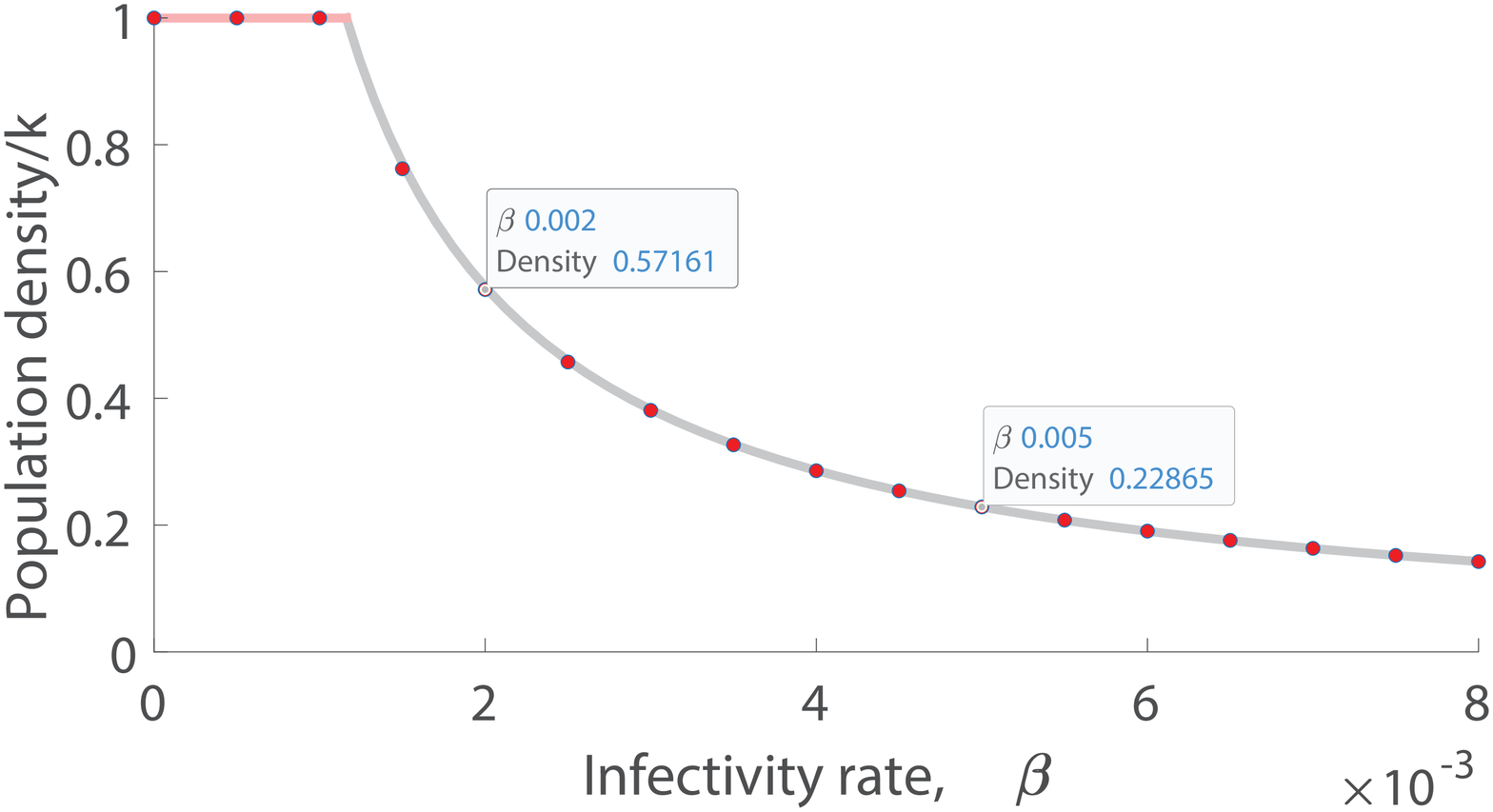}
  \end{subfigure}
  \caption[]{Tumour densities at day end of simulations. (A) At infectivity rate of $\beta = 0.002$, close to branch point, the travelling waves experiences a small period of oscillations and quickly settles into a stable travelling wave. The tail of the wave has a density of 0.57161 which is in agreement with the predicted density in Figure~\ref{fig:Bif_beta}. (B) A plot of the tail densities of the travelling waves, at the end of simulations, for increasing infectivity rate. The wave experiences dampened oscillations which are larger in amplitude and longer in duration before stabilising. The grey line represents the bifurcation analysis in Figure~\ref{fig:Bif_beta}A for $\beta$ values below $\beta_{HB}$. We see that the tail densities from the PDE simulations is agreement with the bifurcation analysis.}
 \label{fig:PDEtoODE1}
\end{figure}

We first compare the results from the one-parameter bifurcation for $\beta$ and tumour population $U$ (Figure~\ref{fig:Bif_beta}) against the spatio-temporal results in Section~\ref{s2} (Figure~\ref{fig:PDEdyn}). The results from the full PDE model demonstrate that for $\beta=0.001$, which is less than the branch point value ($\beta_{BP}=0.00114$), the tumour recovers to full density very quickly and we have a stable travelling wave solution as in Figure~\ref{fig:PDEdyn}A. As predicted by the bifurcation analysis, infection rates below $\beta_{BP}$ will result in full tumour recovery at the stable equilibrium value of $1$. At an increased infectivity rate of 0.002 the tumour population experiences transient oscillations before settling into a stable travelling wave at a reduced density as in Figure~\ref{fig:PDEdyn}B. This result is in agreement with the bifurcation analysis since stable solutions below full density emerge as $\beta$ is increased past the branch point and treatment begins to take effect. For $\beta=0.005$, the transient oscillations in the tumour density increase in amplitude and duration but eventually the travelling wave stabilises at a significantly reduced density. Similarly, the bifurcation diagram shows that stable solutions for tumour density continue to reduce as infectivity is increased. Interestingly, the predicted equilibrium values in tumour density for a given infectivity rate from the bifurcation analysis in Figure~\ref{fig:Bif_beta} can also be observed in the spatial results. In Figure~\ref{fig:PDEtoODE1}A we see that the tail of the travelling wave solutions stabilise at the predicted values from the bifurcation analysis. Specifically, at $\beta = 0.002$, the density in the central region of the tumour is 0.57161. If we plot the tail densities of the travelling wave solutions, for increasing $\beta$ (at the end of simulations) as in Figure~\ref{fig:PDEtoODE1}B, we find that the tumour densities are in agreement with the bifurcation analysis in Figure~\ref{fig:Bif_beta}A, represented by the grey line. This means that our ODE analysis can approximate the impact of virus treatment on tumour outcome by locating the equilibrium densities that stabilise our travelling wave solutions. Note that we extended the domain to $L=80$ and ran simulations up to 500 days to avoid the effects from boundaries and to recover the stable travelling wave after transient oscillations.\\

\begin{figure}[]
  \centering
  \begin{subfigure}[t]{0.02\textwidth}
    \textbf{A}
  \end{subfigure}
  \begin{subfigure}[t]{0.46\textwidth}
    \includegraphics[width=\linewidth, valign=t]{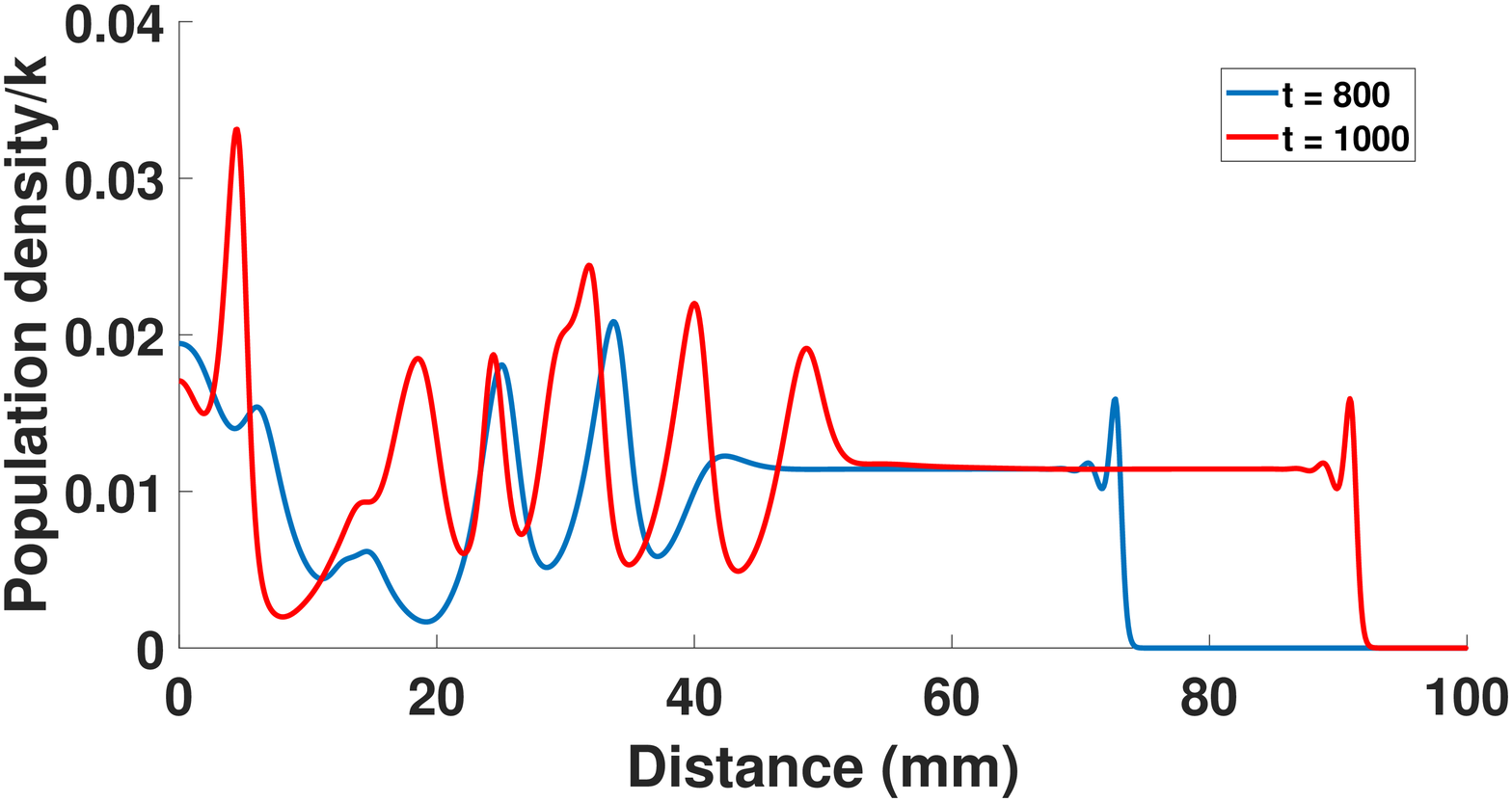}
  \end{subfigure}\hfill
  \begin{subfigure}[t]{0.02\textwidth}
    \textbf{B}
  \end{subfigure}
  \begin{subfigure}[t]{0.46\textwidth}
    \includegraphics[width=\linewidth, valign=t]{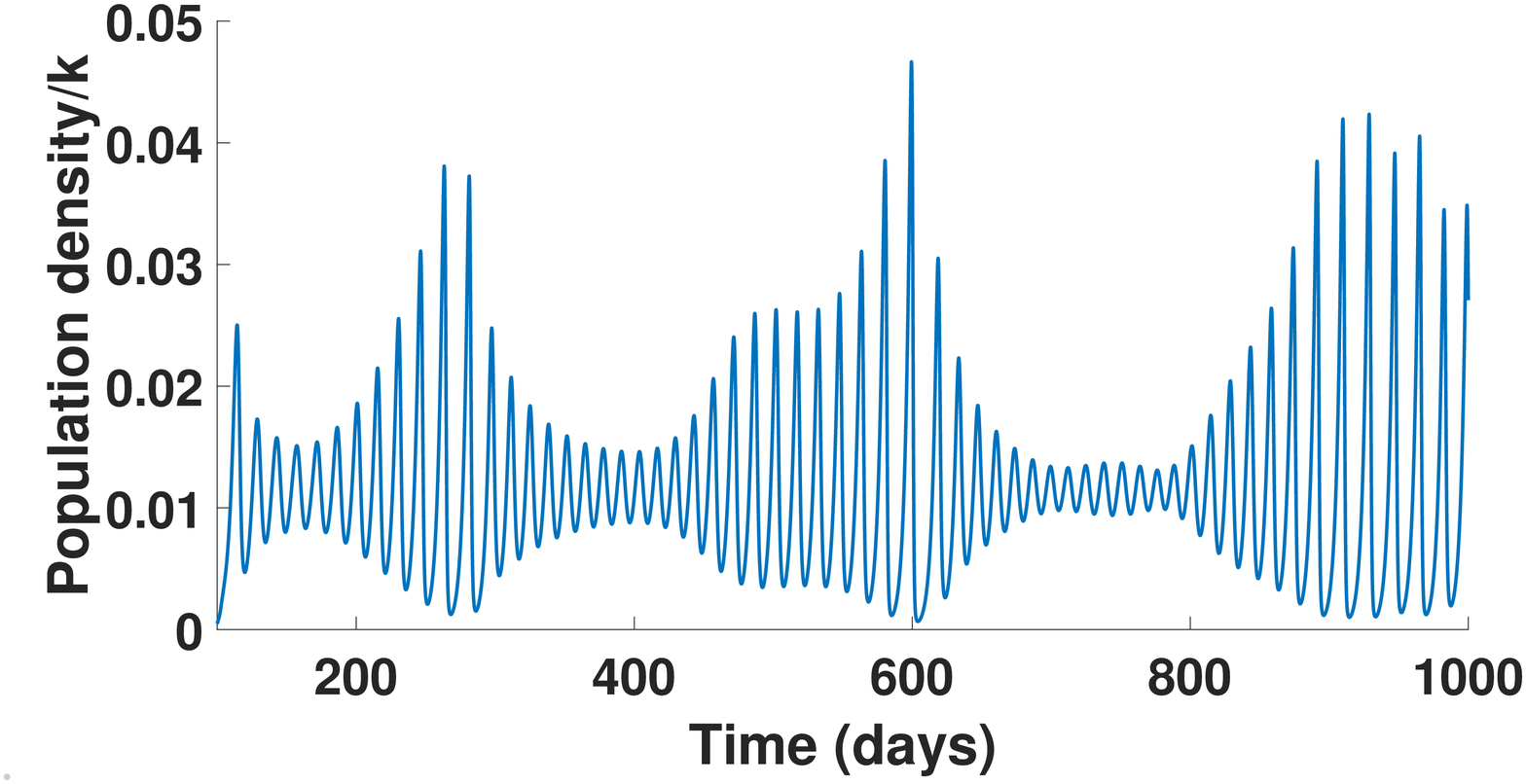}
  \end{subfigure}
  \caption[]{Simulations over 1000 days at $\beta = 0.1 > \beta_{HB}$. (A) Travelling waves have a stable front but continue oscillations in the tail. Choosing a point in space (r = 5 mm) and tracking the oscillations, we see that the oscillations are undampened up to 1000 days. The wave reaches the boundary at approximately 850 days.}
 \label{fig:BetaHB}
\end{figure}

Increasing the infectivity rate even further, the ODE analysis predicts a Hopf bifurcation at $\beta_{HB} = 0.00871$. To see if oscillations in the PDE system persist, we set infectivity at $\beta>\beta_{HB}$ and run the simulation over 1000 days. At $\beta=0.1$, after initial oscillations, the crest of the tumour density wave stabilises; however, the tail of the wave continues to oscillate as in Figure~\ref{fig:BetaHB}A. If we choose a point in space (for example at $r=5$ mm) we can see that the oscillations are not dampening in 1000 days as demonstrated in Figure~\ref{fig:BetaHB}B. Upon reaching the boundary at approximately 850 days, the solution continues to oscillate.\\

\begin{figure}[]
  \centering
  \begin{subfigure}[t]{0.03\textwidth}
    \textbf{A}
  \end{subfigure}
  \begin{subfigure}[t]{0.46\textwidth}
    \includegraphics[width=\linewidth, valign=t]{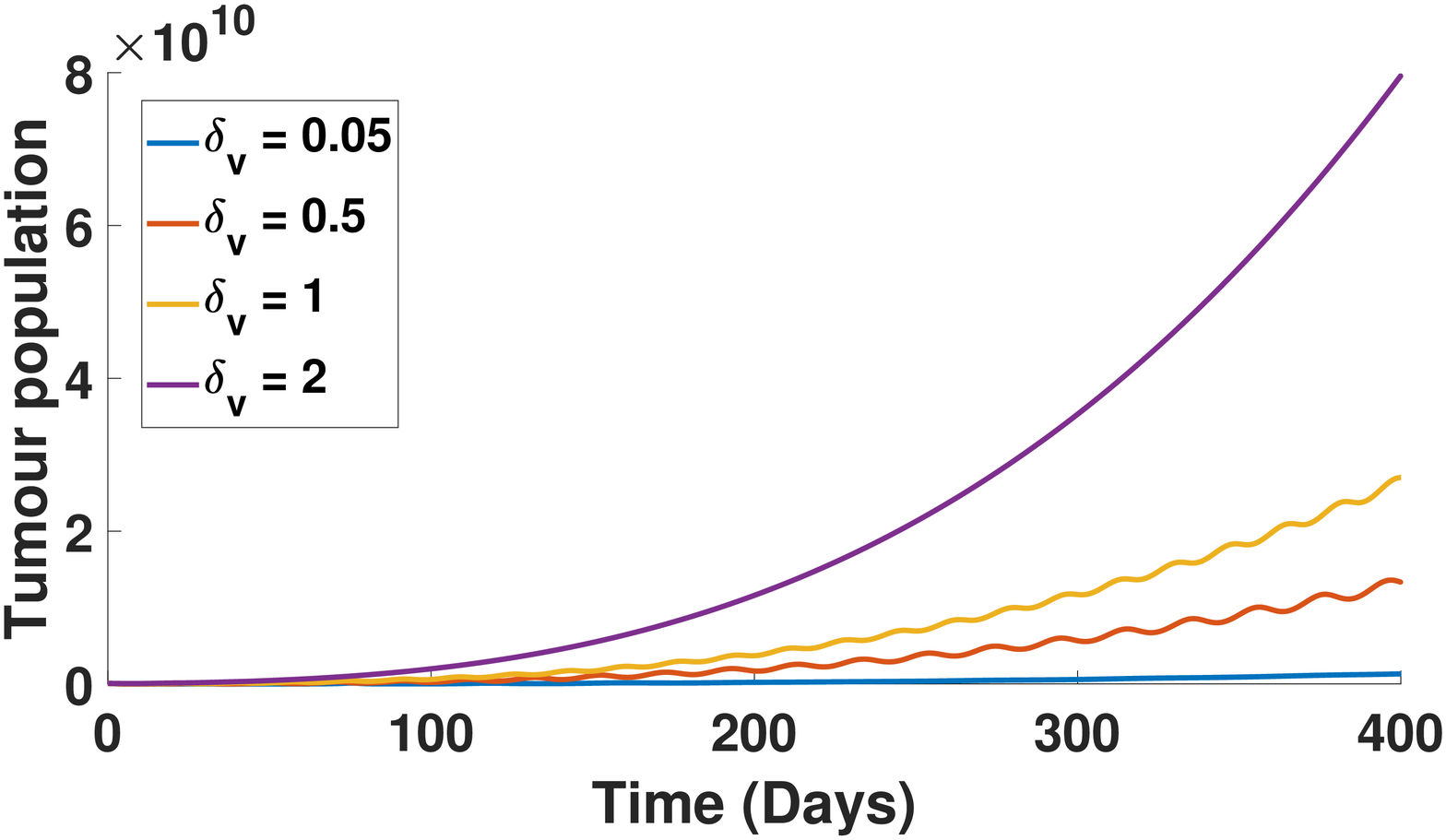}
  \end{subfigure}\hfill
  \begin{subfigure}[t]{0.03\textwidth}
    \textbf{B}
  \end{subfigure}
  \begin{subfigure}[t]{0.46\textwidth}
    \includegraphics[width=\linewidth, valign=t]{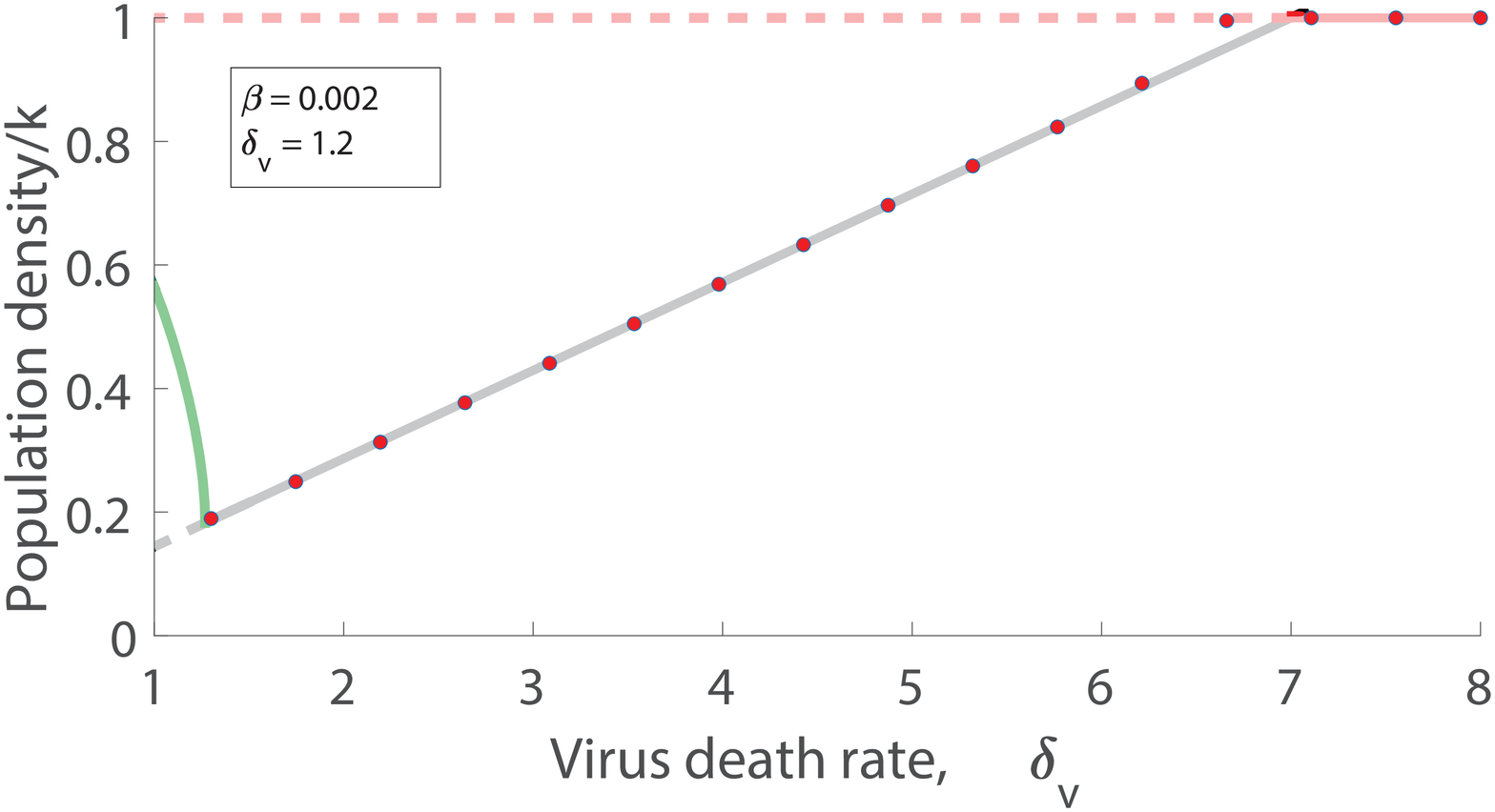}
  \end{subfigure}
  \caption[]{Comparing the spatio-temporal outcomes against the temporal analysis for the virus death rate, $\delta_v$, at infectivity rate, $\beta=0.002$ and infected cell death, $\delta_I = 1.2$. (A) Total tumour population exhibit oscillatory behaviour at virus death rate values which fall in the region of the Hopf bifurcation analysis ($\delta_v \leq 1.276$) in \ref{fig:dr_b}A. (B) Tumour tail densities for virus death rates greater than 1.276. The PDE solutions (dots) reach stable travelling wave solutions, in some cases after transient oscillations. The grey line represents the equilibria from the bifurcation analysis. This plot matches the predicted temporal analysis in Figure~\ref{fig:dr_b}A.}
 \label{fig:ModDeltaI}
\end{figure}

\begin{figure}[]
  \centering
  \begin{subfigure}[t]{0.03\textwidth}
    \textbf{A}
  \end{subfigure}
  \begin{subfigure}[t]{0.46\textwidth}
    \includegraphics[width=\linewidth, valign=t]{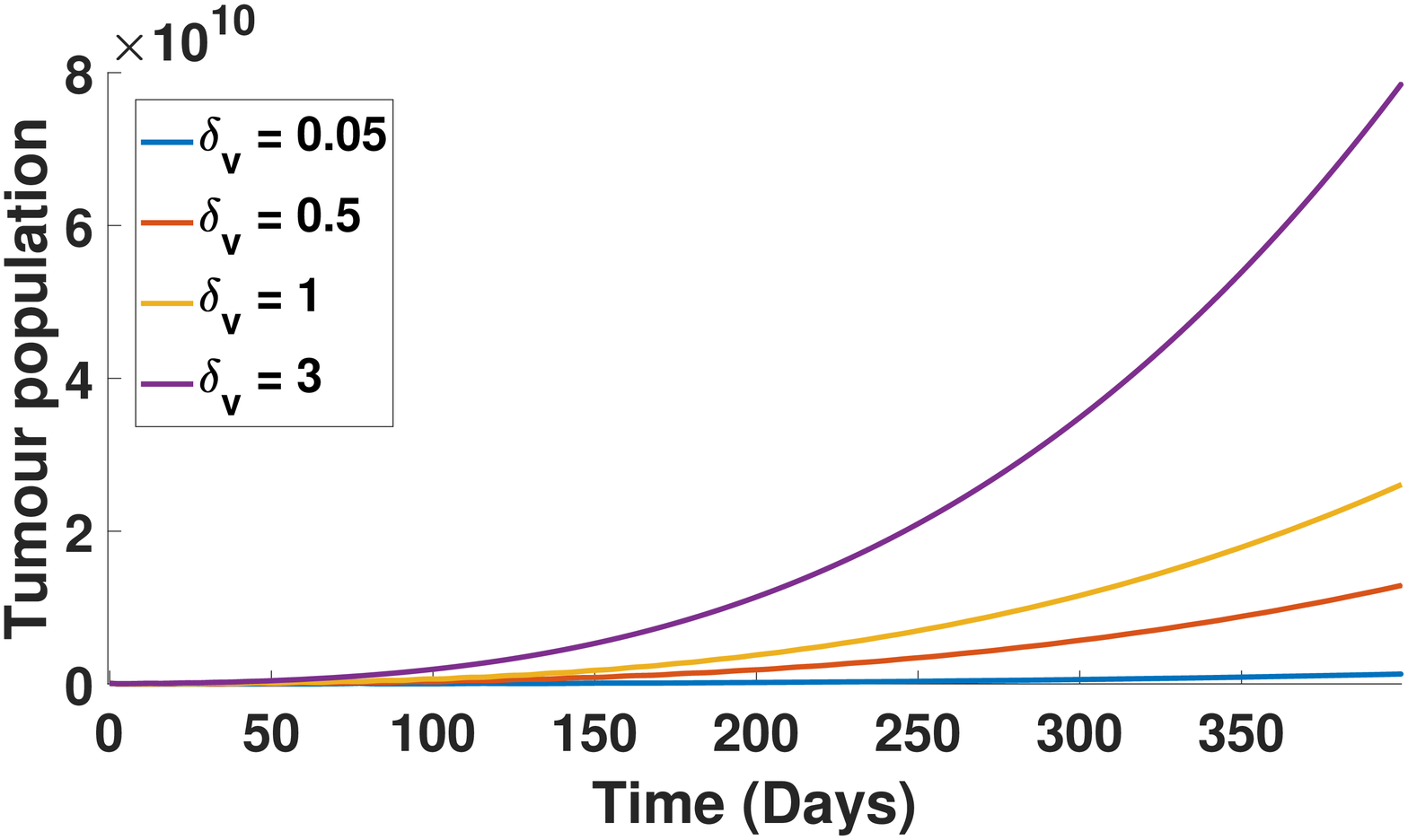}
  \end{subfigure}\hfill
  \begin{subfigure}[t]{0.03\textwidth}
    \textbf{B}
  \end{subfigure}
  \begin{subfigure}[t]{0.46\textwidth}
    \includegraphics[width=\linewidth, valign=t]{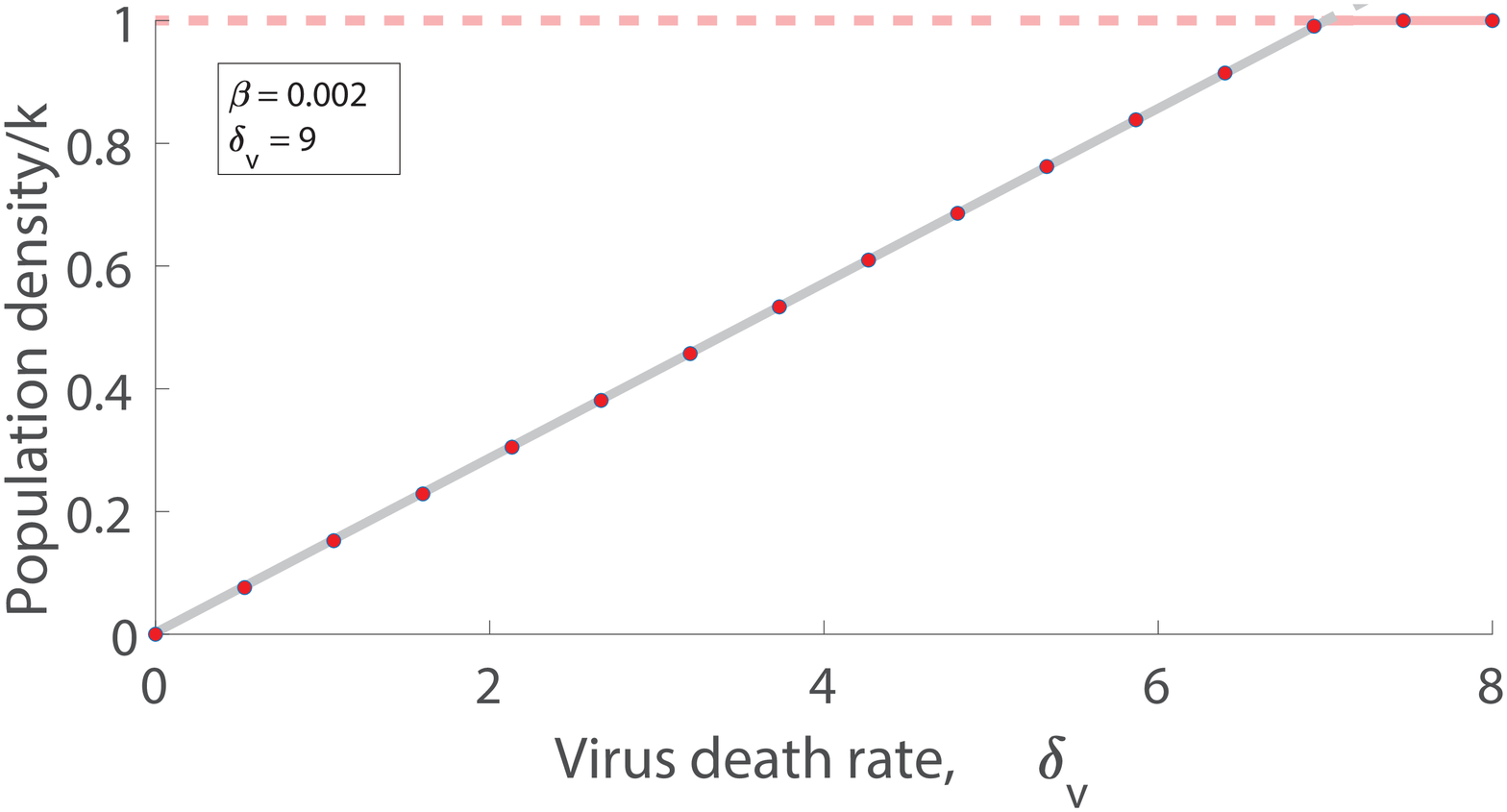}
  \end{subfigure}
  \caption[]{Comparing the spatio-temporal outcomes against the temporal analysis for the virus death rate, $\delta_v$, at infectivity rate, $\beta=0.002$ and high infected cell death, $\delta_I = 9$. (A) There are no oscillation in the total tumour population at low to high virus death rates as predicted by the bifurcation analysis in Figure values which fall in the region of the Hopf bifurcation analysis ($\delta_v \leq 1.276$) in \ref{fig:dr_b}B. (B) Tumour tail densities for $0 \leq \delta_v \leq 8$. The PDE solutions (dots) are stable travelling wave solutions. The grey line represents the equilibria solutions from the bifurcation analysis. This plot matches the predicted temporal analysis in Figure~\ref{fig:dr_b}B.}
 \label{fig:HighDeltaI}
\end{figure}

The analysis in Section~\ref{s3} also investigates changes in virus and infected cell lifespan, $\delta_v$ and $\delta_I$. The dynamics observed in this analysis can also predict the long term behaviour of the PDE system to changes in these parameters. At infectivity rate $\beta=0.002$ and infected cell death rate $\delta_I=1.2$, varying $\delta_v$ in the PDE system produces oscillations in the tumour population when the virus death rate is below the Hopf bifurcation point, $\delta_v \leq 1.276$, refer to Figure~\ref{fig:ModDeltaI}A. For $\delta_v > 1.276$, Figure~\ref{fig:ModDeltaI}B shows that tail densities of the tumour cells (at the end of simulations), are predicted in the bifurcation analysis in Figure~\ref{fig:dr_b}A (grey line in Figure~\ref{fig:ModDeltaI}B). Simulations at a large infected cell death rate $\delta_I = 9$ demonstrate a loss of oscillations regardless of changes in virus death rate, as in Figure~\ref{fig:HighDeltaI}A. The lack of oscillations at a high $\delta_I$ is predicted by the bifurcation analysis in Figure~\ref{fig:dr_b}B with a stable branch connecting two branch points. In a similar manner we can show that the results from the one-parameter analysis of $\delta_v$ at high $\delta_I$ can predict the general behaviour observed in the PDE system as in Figure~\ref{fig:HighDeltaI}B. Overall, the bifurcation analysis can assist in predicting which enhancements will results in the emergence of oscillations and how the general density of the tumour will be affected by the treatment. Stable equilibrium solutions are characterised by stable travelling wave solutions in the full PDE system (at times with transient oscillations) and Hopf bifurcations translate to ongoing oscillations in the tail of the travelling waves.

\section{Discussion}\label{s4}

In this work, a PDE model for virus-tumour dynamics in the context of oncolytic virotherapy has been proposed and analysed. The formulation has spherical symmetry and accounts for virus and cancer diffusion, infection and lysis of tumour cells, and virus decay. Loss of virus particles occurs via infection and internalisation by uninfected and infected tumour cells. 

Particular attention has been placed on infectivity, a model parameter that controls the ability of viral particles to successfully invade and kill tumour cells. This parameter encapsulates different overarching characteristics of the tumour and surrounding cells, such as tumour density, collagen barriers and the expression of receptors on tumour cells to allow virus binding. In general, we observe that for limited increases, higher infectivity leads to better outcomes and more extensive eradication, but full eradication does not seem possible for biologically admissible parameters. Despite an increase in the wave speed of virus propagation through the tumour mass with a larger number of infected tumour cells, no remission is apparent. Once viral characteristic are such that tumour reduction is observed, oscillations in tumour load are always present. These findings are well in agreement with existing experimental evidence, where fluctuations in tumour densities are always emerging after an initial decline due to infection. This occurs in a variety of contexts and for different virus types, suggesting that virotherapy alone may not be sufficient in eradicating solid tumours. A good example is the in vivo experiments from Kim {\em et al.} where no eradication is observed when five different cell lines are treated with both a standard adenovirus and a relaxin-expressing virus, refer to Figure 4 in \cite{kim2006relaxin}.

Nonetheless, benefits can still be achieved in terms of reduction of the existing mass and retardation of proliferating processes. To assess the effect of different viral characteristics on therapy success, an ODE has been obtained from the PDE model by neglecting the diffusive terms and assuming well-mixed populations in the injection zone. The most important feature of this model is that oscillations are widespread and emerge for a large subset of the biologically relevant parameter space. Mathematically, the equations bear similarities with existing models in the literature and orbits show analogous behaviours in terms of shape and period. Uniquely to this model, once limit cycles are born, they persist for all possible increasing values of infectivity and no eradication equilibrium is possible, not even for unrealistically high values. Similarly, a virus that reproduces more abundantly and more frequently in host cells can lead to lower tumour loads, but the reduced ODE model shows that tumour clearance is also not achievable by very high viral burst size, $\alpha$, alone. Further, if the growth is logistic, the average rate of growth only controls the speed at which treatment takes its course, delaying or advancing the final outcome, but does not alter the outcome. This persistence of oscillatory behaviour and absence of tumour-free solutions is mathematically attributable to the role of $\beta$ in Eqs.~(\ref{ode1})-(\ref{ode3}). Infectivity directly controls the loss of virus due to cell infection and any increase in infection is counteracted by a proportional decrease in viral load. This mechanism prevents the model from reaching a stable, full-eradication equilibrium. Note that previous studies mostly concentrate on the behaviour of free virus and do not present a similar term. 

By comparing the results from the PDE system to the ODE analysis, we have find that the dynamics described above hold true in the full PDE model. This allows the use of a simpler temporal model to assist in finding some of the important changes in the behaviour of virus treatment and tumour outcome such as obtaining tumour densities as a result of virus infectivity rate, virus lifespan or burst size.

The model also suggests some interesting considerations that could be useful when a virus for virotherapy is designed and engineered.  Often this is achieved for the purpose of enhancing the features that facilitate infiltration and infection. The model shows that there is a delicate equilibrium among parameters in the system and that care must be used in finding the right balance among those that affect the final outcome. The bifurcation in the ODE system shows, for example, that changes in the rate of death for infected cells can be responsible for the appearance of tumour cell oscillations that are significantly higher than the value of the co-existence equilibrium, where the tumour density is much less than the carrying capacity. Mechanisms aimed at shortening the lifespan of infected cells do not change the co-existence equilibrium but can also result in oscillations that are unwanted, with maxima in tumour densities that increase as lifespan is reduced. Designing therapies where cell death is increased can thus be counterproductive if no other processes in tumour-virus interaction are adjusted. Similarly, the model points to the existence of oscillatory phenomena in cases where the virus lifespan is very long. Although bifurcation analysis indicates that a longer lasting virus can have positive effects in terms of tumour reduction, engineered viruses with highly enhanced duration can also lead to unwanted surges in tumour populations. Within the assumptions of the model, these findings overall suggest that there could exist a tradeoff between the need for enhancing the virus characteristics and tumour ability to respond with a high regrowth rate.

\section{Conclusions}\label{s5}

The model presented in this work investigates the interplay among different aspects of oncolytic virotherapy with a particular focus on those characteristics of the tumour environment and viral loads that influence the success of therapy outcomes. Calibrated using existing experimental data, the model and its analysis shed light on the reasons for the appearance of oscillations in tumour populations when a virus is inoculated in a growing, solid tumour mass. 

The model rests on a number of assumptions that need to be taken into consideration when general conclusions are drawn. Firstly, the use of a phenomenological parameter that captures, on average, the ability of viruses to infect host cells is a limitation. This parameter does not allow for discrimination of finer scale effects related to tumour spatial inhomogeneities or viral penetration. It is also unlikely that internalisation remains constant for the duration of the process, which, in the experiments of Kim {\em et al.}, is around $40$ days. Also, the assumption that the burst rate is unaltered throughout all viral generations is unlikely as virus ability to infect and penetrate cells may degrade. This consideration is also valid for other parameters, such as decay rates or diffusion constants, which are likely to change during the course of the therapy. It is important to stress that this model does not consider the role of the innate immune response (or any immune response). Although this is in line with most in vitro experiments, a model of in vivo dynamics would not be complete without such a response. 

The simple virus diffusion term $D_v$ in this model is also a limitation and cannot describe changes in diffusivity in a heterogeneous tumour environment. Diffusion of viruses may be a function of virus concentration or a function of external factors such as collagen density. All these limitations are currently being addressed and will be part of future work. Nonetheless, this work is able to point out interesting effects that affect the outcome of oncolytic virotherapy, showing a delicate balance among model parameters. Results from this work also seem to suggest that virotherapy alone is unlikely to provide a complete success in tumour eradication and that boosting infectivity cannot be the sole strategy for obtaining a viable protocol for effective anti-tumour responses.

\section{Acknowledgements}
The authors gratefully acknowledge support for this work through the Australian Government Research Training Program Scholarship (PP) and the Australian Research Council Discovery Project DP180101512 (PSK,FF).

\clearpage
\section*{References}

\bibliography{mybibfile}

\begin{thebibliography}{10}
\expandafter\ifx\csname url\endcsname\relax
  \def\url#1{\texttt{#1}}\fi
\expandafter\ifx\csname urlprefix\endcsname\relax\def\urlprefix{URL }\fi
\expandafter\ifx\csname href\endcsname\relax
  \def\href#1#2{#2} \def\path#1{#1}\fi

\bibitem{kelly2007history}
E.~Kelly, S.~J. Russell, History of oncolytic viruses: genesis to genetic
  engineering, Mol. Ther. 15~(4) (2007) 651--659.

\bibitem{patel2013oncolytic}
M.~R. Patel, R.~A. Kratzke, Oncolytic virus therapy for cancer: the first wave
  of translational clinical trials, Transl Res 161~(4) (2013) 355--364.

\bibitem{eissa2018current}
I.~R. Eissa, I.~Bustos-Villalobos, T.~Ichinose, S.~Matsumura, Y.~Naoe,
  N.~Miyajima, D.~Morimoto, N.~Mukoyama, W.~Zhiwen, M.~Tanaka, et~al., The
  current status and future prospects of oncolytic viruses in clinical trials
  against melanoma, glioma, pancreatic, and breast cancers, Cancers 10~(10)
  (2018) 356.

\bibitem{kim2012enhancing}
J.~Kim, P.-H. Kim, S.~W. Kim, C.-O. Yun, Enhancing the therapeutic efficacy of
  adenovirus in combination with biomaterials, Biomaterials 33~(6) (2012)
  1838--1850.

\bibitem{alvarez2012nk}
C.~A. Alvarez-Breckenridge, J.~Yu, R.~Price, J.~Wojton, J.~Pradarelli, H.~Mao,
  M.~Wei, Y.~Wang, S.~He, J.~Hardcastle, et~al., Nk cells impede glioblastoma
  virotherapy through nkp30 and nkp46 natural cytotoxicity receptors, Nat.
  Med.(N.Y.) 18~(12) (2012) 1827--1834.

\bibitem{stohrer2000oncotic}
M.~Stohrer, Y.~Boucher, M.~Stangassinger, R.~K. Jain, Oncotic pressure in solid
  tumors is elevated, Cancer Res. 60~(15) (2000) 4251--4255.

\bibitem{jain2010delivering}
R.~K. Jain, T.~Stylianopoulos, Delivering nanomedicine to solid tumors, Nat.
  Rev. Clin. Oncol. 7~(11) (2010) 653.

\bibitem{mckee2006degradation}
T.~D. McKee, P.~Grandi, W.~Mok, G.~Alexandrakis, N.~Insin, J.~P. Zimmer, M.~G.
  Bawendi, Y.~Boucher, X.~O. Breakefield, R.~K. Jain, Degradation of fibrillar
  collagen in a human melanoma xenograft improves the efficacy of an oncolytic
  herpes simplex virus vector, Cancer Res. 66~(5) (2006) 2509--2513.

\bibitem{kim2006relaxin}
J.-H. Kim, Y.-S. Lee, H.~Kim, J.-H. Huang, A.-R. Yoon, C.-O. Yun, Relaxin
  expression from tumor-targeting adenoviruses and its intratumoral spread,
  apoptosis induction, and efficacy, J. Natl. Cancer Inst. 98~(20) (2006)
  1482--1493.

\bibitem{ganesh2008intratumoral}
S.~Ganesh, M.~Gonzalez-Edick, D.~Gibbons, M.~Van~Roey, K.~Jooss, Intratumoral
  coadministration of hyaluronidase enzyme and oncolytic adenoviruses enhances
  virus potency in metastatic tumor models, Clin. Cancer Res. 14~(12) (2008)
  3933--3941.

\bibitem{guedan2010hyaluronidase}
S.~Guedan, J.~J. Rojas, A.~Gros, E.~Mercade, M.~Cascallo, R.~Alemany,
  Hyaluronidase expression by an oncolytic adenovirus enhances its intratumoral
  spread and suppresses tumor growth, Mol. Ther. 18~(7) (2010) 1275--1283.

\bibitem{diop2011losartan}
B.~Diop-Frimpong, V.~P. Chauhan, S.~Krane, Y.~Boucher, R.~K. Jain, Losartan
  inhibits collagen i synthesis and improves the distribution and efficacy of
  nanotherapeutics in tumors, Proc. Natl. Acad. Sci. U.S.A. 108~(7) (2011)
  2909--2914.

\bibitem{yu2017t}
F.~Yu, B.~Hong, X.-T. Song, et~al., A t-cell engager-armed oncolytic vaccinia
  virus to target the tumor stroma, Cancer Transl. Med. 3~(4) (2017) 122.

\bibitem{de2019targeting}
J.~de~Sostoa, C.~A. Fajardo, R.~Moreno, M.~D. Ramos, M.~Farrera-Sal,
  R.~Alemany, Targeting the tumor stroma with an oncolytic adenovirus secreting
  a fibroblast activation protein-targeted bispecific t-cell engager, J.
  Immunother. Cancer 7~(1) (2019) 1--15.

\bibitem{yoon2015vesicular}
A.-R. Yoon, J.~Hong, C.-O. Yun, A vesicular stomatitis virus glycoprotein
  epitope-incorporated oncolytic adenovirus overcomes car-dependency and shows
  markedly enhanced cancer cell killing and suppression of tumor growth,
  Oncotarget 6~(33) (2015) 34875.

\bibitem{gomez2016combined}
J.~G. Gomez-Gutierrez, J.~Nitz, R.~Sharma, S.~L. Wechman, E.~Riedinger,
  E.~Martinez-Jaramillo, H.~S. Zhou, K.~M. McMasters, Combined therapy of
  oncolytic adenovirus and temozolomide enhances lung cancer virotherapy in
  vitro and in vivo, Virology 487 (2016) 249--259.

\bibitem{meisen2015impact}
W.~H. Meisen, E.~S. Wohleb, A.~C. Jaime-Ramirez, C.~Bolyard, J.~Y. Yoo,
  L.~Russell, J.~Hardcastle, S.~Dubin, K.~Muili, J.~Yu, et~al., The impact of
  macrophage-and microglia-secreted tnf$\alpha$ on oncolytic hsv-1 therapy in
  the glioblastoma tumor microenvironment, Clin. Cancer Res. 21~(14) (2015)
  3274--3285.

\bibitem{russell2019oncolytic}
L.~Russell, K.~W. Peng, S.~J. Russell, R.~M. Diaz, Oncolytic viruses: priming
  time for cancer immunotherapy, BioDrugs (2019) 1--17.

\bibitem{gujar2018antitumor}
S.~Gujar, J.~G. Pol, Y.~Kim, P.~W. Lee, G.~Kroemer, Antitumor benefits of
  antiviral immunity: an underappreciated aspect of oncolytic virotherapies,
  Trends Immunol. 39~(3) (2018) 209--221.

\bibitem{zheng2019oncolytic}
M.~Zheng, J.~Huang, A.~Tong, H.~Yang, Oncolytic viruses for cancer therapy:
  barriers and recent advances, Mol. Ther. Oncolytics 15 (2019) 234--247.

\bibitem{de2020oncolytic}
A.~L. de~Matos, L.~S. Franco, G.~McFadden, Oncolytic viruses and the immune
  system: The dynamic duo, Mol. Ther. Methods Clin. Dev. 17 (2020) 349--358.

\bibitem{hong2019overcoming}
J.~Hong, C.-O. Yun, Overcoming the limitations of locally administered
  oncolytic virotherapy, BMC Biomed. Eng. 1~(1) (2019) 17.

\bibitem{phan2017role}
T.~A. Phan, J.~P. Tian, The role of the innate immune system in oncolytic
  virotherapy, Comput Math Methods Med 2017.

\bibitem{wodarz2004computational}
D.~Wodarz, et~al., Computational approaches to study oncolytic virus therapy:
  insights and challenges, Gene Ther Mol Biol 8 (2004) 137--146.

\bibitem{wodarz2009towards}
D.~Wodarz, N.~Komarova, Towards predictive computational models of oncolytic
  virus therapy: basis for experimental validation and model selection, PLoS
  ONE 4~(1) (2009) e4271.

\bibitem{tian2011replicability}
J.~P. Tian, The replicability of oncolytic virus: defining conditions in tumor
  virotherapy, Math Biosci Eng 8~(3) (2011) 841.

\bibitem{wu2001modeling}
J.~T. Wu, H.~M. Byrne, D.~H. Kirn, L.~M. Wein, Modeling and analysis of a virus
  that replicates selectively in tumor cells, Bull. Math. Biol. 63~(4) (2001)
  731.

\bibitem{wein2003validation}
L.~M. Wein, J.~T. Wu, D.~H. Kirn, Validation and analysis of a mathematical
  model of a replication-competent oncolytic virus for cancer treatment:
  implications for virus design and delivery, Cancer Res. 63~(6) (2003)
  1317--1324.

\bibitem{Fried2006}
A.~Friedman, J.~Tian, G.~Fulci, E.~Chiocca, J.~Wang, Glioma {V}irotherapy:
  {E}ffects of innate immune suppression and increased viral replication
  capacity, Cancer Res. 66 (2006) 2314--2319.

\bibitem{rodriguez2017complex}
I.~A. Rodriguez-Brenes, A.~Hofacre, H.~Fan, D.~Wodarz, Complex dynamics of
  virus spread from low infection multiplicities: implications for the spread
  of oncolytic viruses, PLoS computational biology 13~(1) (2017) e1005241.

\bibitem{jenner2020enhancing}
A.~L. Jenner, F.~Frascoli, A.~C. Coster, P.~S. Kim, Enhancing oncolytic
  virotherapy: Observations from a voronoi cell-based model, J. Theor. Biol.
  485 (2020) 110052.

\bibitem{mok2009mathematical}
W.~Mok, T.~Stylianopoulos, Y.~Boucher, R.~K. Jain, Mathematical modeling of
  herpes simplex virus distribution in solid tumors: implications for cancer
  gene therapy, Clin. Cancer Res. 15~(7) (2009) 2352--2360.

\bibitem{wodarz2012complex}
D.~Wodarz, A.~Hofacre, J.~W. Lau, Z.~Sun, H.~Fan, N.~L. Komarova, Complex
  spatial dynamics of oncolytic viruses in vitro: mathematical and experimental
  approaches, PLoS Comput. Biol. 8~(6) (2012) e1002547.

\bibitem{malinzi2017modelling}
J.~Malinzi, A.~Eladdadi, P.~Sibanda, Modelling the spatiotemporal dynamics of
  chemovirotherapy cancer treatment, J. Biol. Dyn. 11~(1) (2017) 244--274.

\bibitem{friedman2018combination}
A.~Friedman, X.~Lai, Combination therapy for cancer with oncolytic virus and
  checkpoint inhibitor: A mathematical model, PLoS ONE 13~(2) (2018) e0192449.

\bibitem{rajani2015harnessing}
K.~R. Rajani, R.~G. Vile, Harnessing the power of onco-immunotherapy with
  checkpoint inhibitors, Viruses 7~(11) (2015) 5889--5901.

\bibitem{fisher1937}
R.~A. Fisher, The wave of advance of advantageous genes, Ann. Eugen. 7 (1937)
  353--369.

\bibitem{murray2002}
J.~D. Murray, Mathematical {B}iology I: {A}n {I}ntroduction, Springer-Verlag,
  Heidelberg, 2002.

\bibitem{marchini2016overcoming}
A.~Marchini, E.~M. Scott, J.~Rommelaere, Overcoming barriers in oncolytic
  virotherapy with hdac inhibitors and immune checkpoint blockade, Viruses
  8~(1) (2016) 9.

\bibitem{kim2015stem}
J.~Kim, R.~R. Hall, M.~S. Lesniak, A.~U. Ahmed, Stem cell-based cell carrier
  for targeted oncolytic virotherapy: translational opportunity and open
  questions, Viruses 7~(12) (2015) 6200--6217.

\bibitem{roy2013cell}
D.~G. Roy, J.~C. Bell, Cell carriers for oncolytic viruses: current challenges
  and future directions, Oncolytic Virother 2 (2013) 47.

\bibitem{seymour2016oncolytic}
L.~W. Seymour, K.~D. Fisher, Oncolytic viruses: finally delivering, Br. J.
  Cancer 114~(4) (2016) 357--361.

\bibitem{lodish2008molecular}
H.~Lodish, A.~Berk, C.~A. Kaiser, M.~Krieger, M.~P. Scott, A.~Bretscher,
  H.~Ploegh, P.~Matsudaira, et~al., Molecular cell biology, Macmillan, 2008.

\bibitem{cowley2014parallel}
G.~S. Cowley, B.~A. Weir, F.~Vazquez, P.~Tamayo, J.~A. Scott, S.~Rusin,
  A.~East-Seletsky, L.~D. Ali, W.~F. Gerath, S.~E. Pantel, et~al., Parallel
  genome-scale loss of function screens in 216 cancer cell lines for the
  identification of context-specific genetic dependencies, Sci Data 1 (2014)
  140035.

\bibitem{paiva2009multiscale}
L.~R. Paiva, C.~Binny, S.~C. Ferreira, M.~L. Martins, A multiscale mathematical
  model for oncolytic virotherapy, Cancer Res. 69~(3) (2009) 1205--1211.

\bibitem{ganly2000productive}
I.~Ganly, V.~Mautner, A.~Balmain, Productive replication of human adenoviruses
  in mouse epidermal cells, J. Virol. 74~(6) (2000) 2895--2899.

\bibitem{chen2001cv706}
Y.~Chen, T.~DeWeese, J.~Dilley, Y.~Zhang, Y.~Li, N.~Ramesh, J.~Lee,
  R.~Pennathur-Das, J.~Radzyminski, J.~Wypych, et~al., Cv706, a prostate
  cancer-specific adenovirus variant, in combination with radiotherapy produces
  synergistic antitumor efficacy without increasing toxicity, Cancer Res.
  61~(14) (2001) 5453--5460.

\bibitem{wodarz2001viruses}
D.~Wodarz, Viruses as antitumor weapons, Cancer Res. 61~(8) (2001) 3501--3507.

\bibitem{bajzer2008modeling}
{\v{Z}}.~Bajzer, T.~Carr, K.~Josi{\'c}, S.~J. Russell, D.~Dingli, Modeling of
  cancer virotherapy with recombinant measles viruses, J. Theor. Biol. 252~(1)
  (2008) 109--122.

\bibitem{komarova2010ode}
N.~L. Komarova, D.~Wodarz, {ODE} models for oncolytic virus dynamics, J. Theor.
  Biol. 263~(4) (2010) 530--543.

\bibitem{dingli2009dynamics}
D.~Dingli, C.~Offord, R.~Myers, K.-W. Peng, T.~W. Carr, K.~Josic, S.~J.
  Russell, Z.~Bajzer, Dynamics of multiple myeloma tumor therapy with a
  recombinant measles virus, Cancer Gene Ther. 16~(12) (2009) 873--882.

\bibitem{wodarz2003gene}
D.~Wodarz, Gene therapy for killing p53-negative cancer cells: use of
  replicating versus nonreplicating agents, Hum. Gene Ther. 14~(2) (2003)
  153--159.

\bibitem{jenner2018minimal}
A.~L. Jenner, A.~C.~F. Coster, P.~S. Kim, F.~Frascoli, Treating cancerous cells
  with viruses: insights from a minimal model for oncolytic virotherapy, Lett
  Biomath 5~(sup1) (2018) S117--S136.

\bibitem{jenner2019gomp}
A.~L. Jenner, P.~S. Kim, F.~Frascoli, Oncolytic virotherapy for tumours
  following a gompertz growth law, J. Theor. Biol. 480 (2019) 129—140.

\bibitem{shashkova2009characterization}
E.~V. Shashkova, S.~M. May, M.~A. Barry, Characterization of human adenovirus
  serotypes 5, 6, 11, and 35 as anticancer agents, Virology (Lond) 394~(2)
  (2009) 311--320.

\end{thebibliography}

\newpage
\appendix
\section{Parameter estimates}
\label{appendix:parameterestimates}

To compare our results with available exponential data, a calibration of the model is necessary. Firstly, we consider that in a solid tumour of radius $1$ mm there are approximately $10^6$ cells of epithelial origin \cite{lodish2008molecular}, so that the carrying capacity of the model is set to $k = 10^6$ cells/mm$^3$. The largest admissible tumour size $L$ is estimated by observing that the survival rate of mice with established U343 tumours, treated with PBS solution in Ref.~\cite{kim2006relaxin}, drops to zero at a tumour volume of 2500 mm$^3$, which corresponds to a radius of 8.4 mm (refer to Figure A U343 in~\cite{kim2006relaxin}). We use a maximum radius of 10 mm in our simulations.
In \cite{kim2006relaxin}, the start of treatment, i.e. inoculation of viral load, is set at a tumour size of $70$ mm$^3$, which, given the spherical symmetry, corresponds to an initial radius of $2.6$ mm. At inoculation, the tumour is assumed to be homogeneous and at carrying capacity, so that the total (uninfected) cell number is $70 \times 10^6$, given a density $U_0 = k$. To estimate the tumour growth rate $r_U$ and diffusion rate $D_u$, the glioblastoma U343 cell line data in Fig.~4A of \cite{kim2006relaxin} is used. The doubling time is approximately $45$ hours ($1.875$ days)~\cite{cowley2014parallel}, so that $r_U = \ln{2}/1.875 \approx 0.3$/day. 

To determine $D_u$ we note that the tumour volume on day $40$ for the U343 control group (PBS) is approximately $1000$ mm$^3$ and employ the formula for volume used again in \cite{kim2006relaxin}, i.e. $\text{Vol} = 0.523L W^2$, where $L=W=2r$. This implies that the radius grows from an initial $2.6$ mm  to approximately $6$ mm in $40$ days. The control group is an untreated tumour, so Eq.~(\ref{eq:TumourPDE}) can be simplified, in the absence of viral loads, to 
\[
\frac{\partial U}{\partial t}  = \frac{D_u}{r^2} \frac{\partial}{\partial r} \Big(r^2 \frac{\partial U}{\partial r} \Big) + r_U U \Big( 1- \frac{U+I}{k} \Big),
\] which is characterised by a wave speed of $c = 2\sqrt{r_U D_u}$ \cite{murray2002}. This allows us to calculate $c$ as the change in radial distance over $40$ days as $c = (6-2.6)/40 = 0.085$ mm/day, and finally gives $D_u = \displaystyle \left(\frac{0.085}{2}\right)^2\times \frac{1}{0.3} \approx 0.006$ mm$^2$/day. 
\begin{figure}
    \centering
    \includegraphics[width=\linewidth]{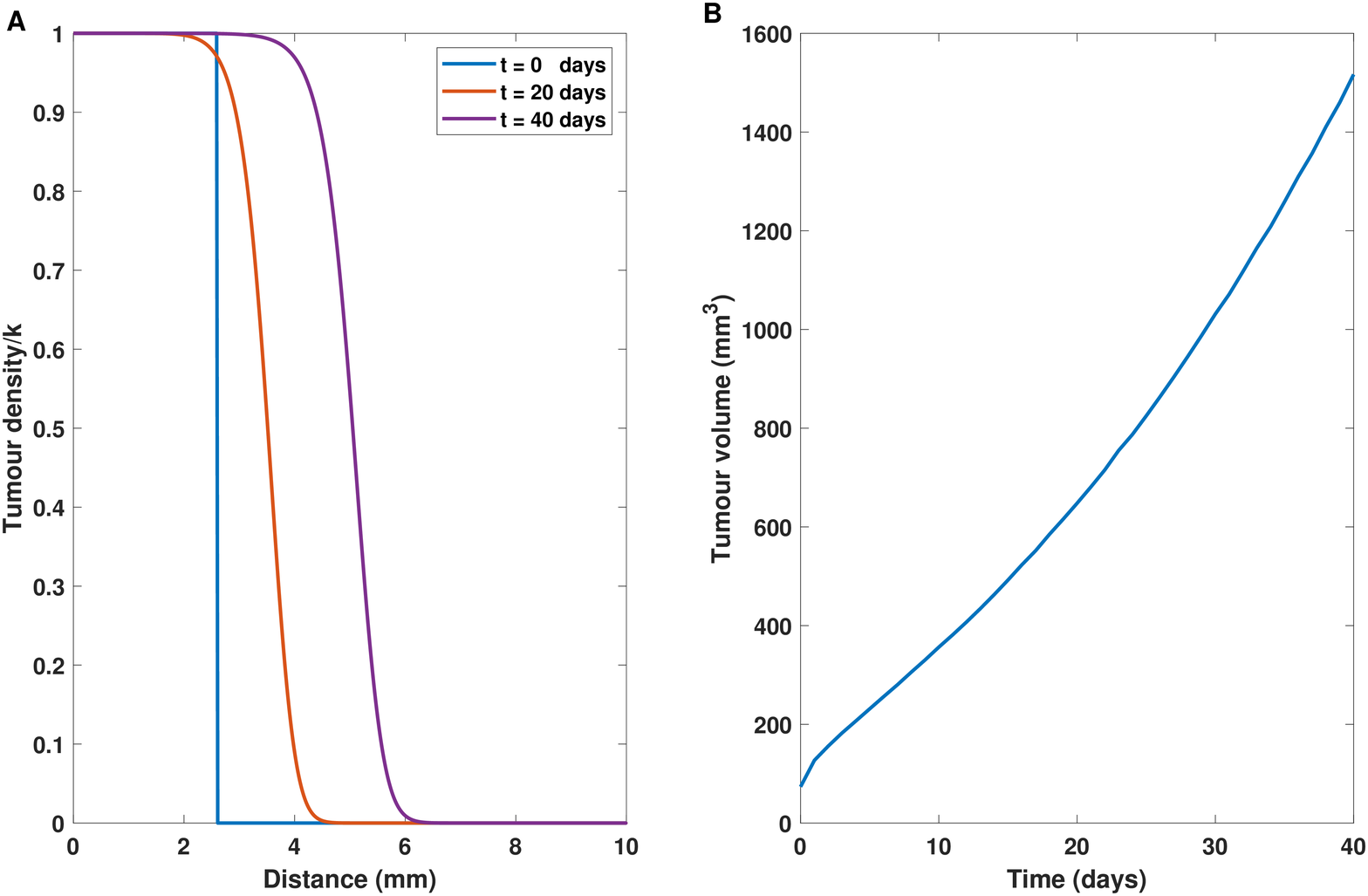}
    \caption{The case of an untreated tumour. (A): Tumour density rescaled by carrying capacity $k$ and (B) tumour volume over 40 days. Note that the tumour grows with a standard travelling wave front.}
    \label{fig:Untreated tumour growth}
\end{figure}
Depending of the virus considered, viral loads can vary. In the experiments by Kim {\em et al.}, the viral injection doses contain approximately $10^{10}$ viral particles. So, assuming that the injections are in a spherical region with radius $0.5$ mm, we use $V_0 = 10^{10}/(\frac{4}{3}(0.5^3) \pi) = 1.9 \times 10^{10}$ viruses/mm$^3$. Ref.~\cite{mok2009mathematical} estimates the degradation of viruses at approximately $4.8~\times 10^{-5}$/s, so that $\delta_v = 4.1$/day. For viral diffusivity we use the same value employed in \cite{paiva2009multiscale}, i.e. $D_v = 0.24$ mm$^2$/day. Infected cells undergo lysis approximately $24$ hours after infection \cite{ganly2000productive}, so $\delta_I = 1$/day. Shashkova {\em et al.}~\cite{shashkova2009characterization} analyse the burst size of various oncolytic adenoviruses in cancer cell lines and conclude that infectious units (IU) per cell is in the range of $1,000$ to $100,000$ IU/cell. Therefore, we use Chen {\em et al.}'s measured virus burst size of 3,500 viral particles per cell~\cite{chen2001cv706}, setting $\alpha = 3500$.

To estimate a biologically acceptable value for $\beta$, data from \cite{Fried2006} are used, where a load of $10^{-9}$ mm$^3$ of virus is shown to infect approximately $70\%$ of available cells in one hour. Assuming no significant degradation occurs in one day except loss due to viral death, and for viral load densities considered in the model, we obtain a value of $\beta \approx 1.5$ x $10^{-9}$ mm$^3$/(viruses $\times$ day). 

Solving Eqs.~(\ref{eq:TumourPDE}-\ref{eq:InfectedPDE}) numerically, the absence of initial treatment, i.e. $V_0 = 0$,  causes the tumour to grow and expand as a standard travelling wave as shown in Fig.~\ref{fig:Untreated tumour growth}(A). The evolution of the volume of the cancer mass for $40$ days (see Fig.~\ref{fig:Untreated tumour growth}(B)) very closely follows the experimental results from Fig.~4A of Ref.~\cite{kim2006relaxin} with similar values and growth pattern. In particular, note the change in velocity of growth after approximately the first $2$ days that matches experimental observation well.

\section{Parameter-dependent equilibria and the case $\delta_V = 0$}
\label{appendix:parameterdependentequilibria}

To determine the equilibria for Eqs.~(\ref{ode1})--(\ref{ode3}) , we set $ \dot{U} = \dot{V} = \dot{I} = 0 $ and solve the resulting equations simultaneously to arrive at the following system:
\begin{align}
    U & \left(r_U (\delta_I + \beta V) U + \delta_I (\beta V - r_U) \right) = 0, \label{simul1}\\
    V & \left( -\beta^2UV + \alpha \beta \delta_I U - \delta_I \delta_v - \beta \delta_I U \right) = 0. \label{simul2}
\end{align}
Determining $V$ as a function of $U$ in the first equation, we can substitute in the second and arrive at the quadratic expression $AU^2 + BU + C=0,$ with 

\begin{align*}
A & = \alpha \beta^2 r \delta_I,\\
B & = \delta_I \beta \left( \alpha \beta \delta_I - r \delta_V - \beta \delta_I - r \beta \right),\\
C & = - \delta_I^2 \delta_V \beta.
\end{align*}
One root corresponds to the $U$-component of solution ${\bf U}_s$ discussed in the main text and represents the density of tumour at the co-existing equilibrium. The other root always assumes negative values for all biologically relevant parameters.  Substituting these solutions into one of Eqs.~(\ref{simul2}), we obtain values for $V$. Using, for example, Eq.~(\ref{ode3}) for $\dot{I}=0$ we arrive at the final expression for ${\bf U}_s$:
\begin{align*}
U_{s} &=\frac{(1-\alpha ) \beta  \delta_I + r (\beta +\delta_V)+
\sqrt{(\alpha -1)^2
   \beta ^2 \delta_I^2+r^2 (\beta +\delta_V)^2+2 \beta  \delta_I r ((1-\alpha ) \beta +(\alpha +1) \delta_V)}}{2 \alpha  \beta  r}, \\
V_{s} &=\frac{(\alpha-1) \beta  \delta_I+r (\beta +\delta_V)-
\sqrt{(\alpha -1)^2
   \beta ^2 \delta_I^2+r^2 (\beta +\delta_V)^2+2 \beta  \delta_I r ((1-\alpha ) \beta +(\alpha +1) \delta_V)}}{2  \beta^2}, \\
I_{s} &=\frac{(\alpha -1) (\beta  r+1)-(\alpha +1) \delta_V r }{2 \alpha  \beta  r} +\\
& \quad +\frac{(\alpha -1)\left(\beta\delta_I(1-\alpha) +
\sqrt{(\alpha -1)^2
   \beta ^2 \delta_I^2+r^2 (\beta +\delta_V)^2+2 \beta  \delta_I r ((1-\alpha ) \beta +(\alpha +1) \delta_V)}\right)}{2 \alpha  \beta  r}.
\end{align*}
A similar expression involving different signs among the terms exists for the other nonbiological equilibrium.

Finally, let us discuss the behaviour of the full-eradication equilibrium $(U,V,I) = (0,0,0)$ for the limiting case $\delta_V=0$. Eqs.~(\ref{ode1})-(\ref{ode3}) change to
\begin{align*}
\dot{U} & =  r_U U \left( 1- (U+I)\right) - \beta UV, \\
\dot{V} & =  \alpha \delta_I I  - \beta (U + I) V, \\
\dot{I} & =  \beta UV - \delta_I I,
\end{align*}
and the original equilibrium $(0,0,0)$ changes to a new equilibrium $(0,V,0)$, which is valid for any initial (admissible) value of the density of infected cells as long as $U=I=0$. The eigenvalues associated to this new solution for the Jacobian from Eq.~(\ref{JacobianODE}) with $\delta_V=0$ are given by $(0, r - V\beta, -\delta_i)$. Hence, any solution for which $V > \beta/r$ yields two negative and one zero eigenvalues. Note that this solution corresponds to the case of an immortal virus, which grows up to the density needed to destroy the whole tumour and then lives forever. 
\end{document}